\documentclass[12 pt]{amsart}
\usepackage{a4}
\usepackage{amssymb}
\DeclareMathOperator{\RE}{Re} 
\DeclareMathOperator{\IM}{Im}

\newtheorem{thm}{Theorem}[section]
\newtheorem{lem}[thm]{Lemma}
\newtheorem{prop}[thm]{Proposition}
\newtheorem{cor}[thm]{Corollary}
\theoremstyle{definition}
\newtheorem{defn}[thm]{Definition}
\newtheorem{example}[thm]{Example}
\theoremstyle{remark}
\newtheorem{remark}[thm]{Remark}

\numberwithin{equation}{section}

\begin{document}

\title[Bishop's theorem and Differentiability]{Bishop's Theorem and
Differentiability of a subspace of $C_b(K)$}

\author{Yun Sung Choi}
\address{Department of Mathematics, POSTECH, San 31, Hyoja-dong,
Nam-gu, Pohang-shi, Kyungbuk, Republic of Korea, +82-054-279-2712}

\email{mathchoi@postech.ac.kr}
\thanks{This work was supported by grant  No. R01-2004-000-10055-0 from the
Basic Research Program of the Korea Science \& Engineering
Foundation}

\author{Han Ju Lee}
\email{hahnju@postech.ac.kr}

\author{Hyun Gwi Song}
\email{hyuns@postech.ac.kr}








\begin{abstract}
Let $K$ be a Hausdorff space and $C_b(K)$ be the Banach algebra of
all complex bounded continuous functions on $K$. We study the
G\^{a}teaux and Fr\'echet differentiability of subspaces of
$C_b(K)$. Using this, we show that the set of all strong peak
functions in a nontrivial separating separable subspace $H$ of
$C_b(K)$ is a dense $G_\delta$ subset of $H$, if $K$ is compact.
This gives a generalized Bishop's theorem, which says that the
closure of the set of strong peak point for $H$ is the smallest
closed norming subset of $H$. The classical Bishop's theorem was
proved for a separating subalgebra $H$ and a metrizable compact
space $K$.

In the case that $X$ is a complex Banach space with the
Radon-Nikod\'ym property, we show that the set of all strong peak
functions in $A_b(B_X)=\{ f\in C_b(B_X) : f|_{B_X^\circ} \mbox{ is
holomorphic}\}$ is dense. As an application, we show that the
smallest closed norming subset of $A_b(B_X)$ is the closure of the
set of all strong peak points for $A_b(B_X)$. This implies that the
norm of $A_b(B_X)$ is G\^{a}teaux differentiable on a dense subset
of $A_b(B_X)$, even though the norm is nowhere Fr\'echet
differentiable when $X$ is nontrivial. We also study the denseness
of norm attaining holomorphic functions and polynomials. Finally we
investigate the existence of numerical Shilov boundary.
\end{abstract}

\maketitle

\section{Introduction}
Let $K$ be a Hausdorff topological space. A {\it function algebra}
$A$ on $K$ will be understood to be a closed subalgebra of $C_b(K)$
which is the Banach algebra of all bounded complex-valued continuous
functions on $K$. The norm $\|f\|$ of a bounded continuous function
$f$ on $K$  is defined to be  $\sup_{x\in K} |f(x)|$. A function
algebra $A$ is called {\it separating} if for two distinct points
$s, t$ in $K$, there is $f\in A$ such that $f(s)\neq f(t).$

In this paper, a subspace means a closed linear subspace. For each
$t\in K$, let $\delta_t$ be an evaluation functional on $C_b(K)$,
that is, $\delta_t(f) = f(t)$ for every $f\in C_b(K)$. A subspace
$A$ of $C_b(K)$ is called {\it separating} if for distinct points
$t, s$ in $K$ we have $\alpha \delta_t \neq \beta \delta_s$ for any
complex numbers $\alpha, \beta$ of modulus $1$ as a linear
functional on $A$. This definition of a separating subspace is a
natural extension of the definition of a separating function
algebra. In fact, given a separating function algebra $A$ on $K$ and
given two distinct points $t, s$ in $K$, we have $\alpha \delta_t
\neq \beta \delta_s$ on $A$ for any nonzero complex numbers $\alpha$
and $\beta$. Otherwise, there are some nonzero complex numbers
$\alpha$ and $\beta$ such that $\alpha \delta_t = \beta \delta_s$ on
$A$. Let $\gamma = \beta/\alpha$. Choose $f\in A$ so that $f(s)=1$
and $f(t)\neq 1$. By assumption, $f(t)= \gamma$. Fix a positive
number $r$ with $0<r<1/\|f\|$. Taking
\[ g= \frac{1-r}{r}\sum_{m=1}^\infty r^mf^m = \frac{(1-r)f}{1-rf},\]
we have $g\in A$ and $g(s)=1$, which imply $g(t)=\gamma=f(t)$. Hence
\begin{equation}\label{eq:sameq1}\gamma=g(t) =
\frac{(1-r)f(t)}{1-rf(t)} = \frac {(1-r)\gamma}{1-r\gamma}.
\nonumber\end{equation} This equation shows that $\gamma =1$ and
$f(t)=1$, which is a contradiction.

A nonzero element $f\in A$ is called a {\it peak function} if there
exists only one point $x\in K$ such that $|f(x)|=\|f\|$. In this
case, the corresponding point $x$ is said to be a {\it peak point}
for $A$. A nonzero element $f\in A$ is called a {\it strong peak
function} if there exists only one point $x\in K$ such that
$|f(x)|=\|f\|$ and for any neighborhood $V$ of $x$, there is
$\delta>0$ such that if $y \in K\setminus V$,  then
 $|f(y)|\le\|f\|-\delta$. In this case, the corresponding point $x$
is called a {\it  strong peak point} for $A$. We denote by $\rho A$
the set of all strong peak point for $A$. Note also that if $K$ is a
compact Hausdorff space, every peak function (resp. peak point for
$A$) is a strong peak function (resp. strong peak point for $A$).

A subset $F$ of $K$ is said to be a {\it norming} subset for $A$ if
for every $f\in A$, we have
\[ \|f\| = \sup_{x\in F} |f(x)|.\] Note that every closed norming
subset contains all strong peak points. If $K$ is a compact
Hausdorff space, then a closed subset $T$ of $K$ is a norming subset
for $A$ if and only if $T$ is a boundary for $A$, that is, for every
$f\in A$, we have
\[ \max_{t\in T} |f(t)| = \|f\|.\]

A famous theorem of Shilov (see \cite{GRS, L}) asserts that if $A$
is a separating function algebra $A$ on a compact Hausdorff space
$K$, then there is a smallest closed boundary for $A$, which is
called the Shilov boundary for $A$ and denoted by $\partial A$. We
shall say that a subspace $A$ of $C_b(K)$ on a Hausdorff space $K$
has the {\it Shilov boundary} if there is a smallest closed norming
subset for $A$. If $K$ is not compact, a separating function algebra
$A$ on $K$ need not have the Shilov boundary (see \cite{ACLQ, CH,
Gl, GM}).

Let $X$ be a real or complex Banach space. We denote by $B_X$ and
$S_X$ the closed unit ball and unit sphere of $X$, respectively. The
norm $\|\cdot \|$ of $X$ is said to be {\it G\^{a}teaux
differentiable} (resp. Fr\'echet differentiable) at $x$ if
\[ \lim_{t\to 0} \frac{\|x+ ty\| + \|x-ty\|-2\|x\|}{t}=0\]  for every
$y\in X$ (resp. uniformly for $y\in S_X$). Notice that if the norm
of a nontrivial Banach space is G\^{a}teaux (resp. Fr\'echet)
differentiable at $x$, then $x\neq 0$ and the norm is also
G\^{a}teaux (resp. Fr\'echet) differentiable at $\alpha x$ for any
nonzero scalar $\alpha$ (For more details, see \cite{DGZ}).

Let $C$ be a convex subset of a complex Banach space. An element
$x\in C$ is said to be an (resp. complex) extreme point of $C$ if
for every $y\neq 0$ in $X$, there is a real (resp. complex) number
$\alpha$, $|\alpha|\le 1$ such that $x+\alpha y\not\in C$. The set
of all (resp. complex) extreme points of $C$ is denoted by $ext(C)$
(resp. $ext_\mathbb{C}(C)$). A point $x^*\in B_{X^*}$ is said to be
a {\it weak-$*$ exposed point} of $B_{X^*}$ if there exists $x\in
S_X$ such that
\[ 1=\RE x^*(x) > \RE y^*(x), \ \ \ \ \forall y^*\in
B_{X^*}.\] The corresponding point $x\in S_X$ is said to be  a {\it
smooth point} of $B_X$. It is well-known \cite{DGZ} that the norm of
$X$ is G\^ateaux differentiable at $x\in S_X$ if and only if $x$ is
a smooth point of $B_X$. The set of all weak-$*$ exposed points of
$B_{X^*}$ is denoted by $w^{\!*\!}exp (B_{X^*})$. It is easy to
check that $w^{\!*\!}exp (B_{X^*})\subset ext (B_{X^*})$.

We denote by $C_b(K, Y)$ the Banach space of all bounded continuous
functions of a Hausdorff space $K$ into a Banach space $Y$ with the
sup norm. By replacing the absolute value with the norm of $Y$, the
notion of a strong peak function or a norming subset for $C_b(K,Y)$
can be defined in the same way as for $C_b(K)$. Note that $\rho
C_b(K)=\rho C_b(K,Y)$ for a nontrivial Banach space $Y$.

Given complex Banach spaces $X, Y$, the following two subspaces of
$C_b(B_X,Y)$ are the ones of our main interest:
\begin{align*}A_b(B_X, Y)&= \{ f\in C_b(B_X, Y): f \text{ is analytic
on the interior of $B_X$} \}\\
A_u(B_X, Y) &= \{ f\in A_b(B_X, Y): f\text{ is uniformly continuous
on } B_X\}.\end{align*} It was shown in \cite{ACG} that these two
function Banach spaces are the same if and only if $X$ is finite
dimensional. Hereafter, $A(B_X, Y)$ will represent either $A_u(B_X,
Y)$ or $A_b(B_X, Y)$, and we simply write $A(B_X)$ instead of
$A(B_X, \mathbb{C})$. For the basic properties of holomorphic
functions on a Banach space, see \cite{ACG, CCG, C, D}. Note that
$\rho A(B_X) = \rho A(B_X, Y)$ for a nontrivial complex Banach space
$Y$.

In Section 2, we find a necessary and sufficient condition of $f$ in
a subspace $A$ of $C_b(K)$ under which the norm is either
G\^{a}teaux or Fr\'{e}chet differentiable at $f$. The main result of
this section is that if $f$ is a strong peak function in $A$, then
the norm of $A$ is G\^{a}teaux differentiable at $f$, and the
converse is also true for a nontrivial separating subspace $A$ of
$C(K)$ on a compact Hausdorff space $K$. Applying them to $A(B_X)$,
we show that the norm of $A(B_X)$ is nowhere Fr\'echet
differentiable, if $X$ is nontrivial. The relation between a
norm-attaining $m$-homogeneous polynomial and its differentiability
 was studied in \cite{F}.

In Section 3, we give another version of Bishop's theorem. If $A$ is
a nontrivial separating separable subspace of $C(K)$ on a compact
Hausdorff space $K$, then the set of all strong peak functions is a
dense $G_\delta$-subset of $A$. Using this fact, we obtain Bishop's
theorem which says that if $A$ is a nontrivial separating separable
subspace of $C(K)$, then $\rho A$ is a norming subset for $A$ and
its closure is the Shilov boundary for $A$.

Globevnik \cite{Gl} studied  norming subsets for $A(B_X)$, when
$X=c_0$. In that paper he showed that neither $A_u(B_X)$ nor
$A_b(B_X)$ has the Shilov boundary. In \cite{ACLQ}, it was shown
that $\partial A(B_X)=S_{X}$ for $X=\ell_p$, $1\le p<\infty$. This
result was generalized in \cite{CHL} to show that $\partial
A(B_X)=S_X$ for a locally uniformly $c$-convex, order continuous
sequence space $X$ (For more details on the $c$-convexity and order
continuity of a Banach lattice, see \cite{CHL, DGT, DHM, Lee2, Lee,
Lee3}).

In Section 4, it is shown that if $X$ has the Radon-Nikod\'ym
property and $Y$ is a nontrivial complex Banach space, then the set
of all strong peak functions in $A(B_X, Y)$ is dense in $A(B_X,Y)$.
Applying this fact, it is also proved that if $X$ has the
Radon-Nikod\'ym property and $Y$ is nontrivial, then $\rho A(B_X)$
is a norming subset for $A(B_X,Y)$. In particular, $\partial
A(B_X,Y)$ is the closure of $\rho A(B_X)$, and
$ext_{\mathbb{C}}(B_X)$ is also a norming subset for $A(B_X,Y)$.

It is worth-while to remark that Bourgain-Stegall's perturbed
optimization principle \cite{S} is the key method to prove these
facts. This method has been used to study the density of the
norm-attaining $m$-homogeneous polynomials and holomorphic functions
on $X$, when  $X$ has the Radon-Nikod\'ym property (see \cite{AAGM,
Bo, CK}).

In Section 5, we modify the argument of Lindenstrauss \cite{Li} with
strong peak points and also with uniformly strongly exposed points,
and show the density of norm-attaining elements in certain subspaces
of $C_b(K,Y)$. We also extend the result of \cite{ArGM} to the
vector valued case by changing their proof, which is based on that
of Lindenstrauss.

In the last section, we apply Bishop's theorem to study  numerical
boundaries for subspaces of $C_b(B_X, X)$. The notion of a numerical
boundary was introduced and studied for various Banach spaces $X$ in
\cite{AK}, and it was observed that the smallest closed numerical
boundary, called the {\it numerical Shilov boundary}, doesn't exist
for some Banach spaces.  We show that there exist the numerical
Shilov boundaries for most subspaces of $C_b(B_X, X)$, if $X$ is
finite dimensional, which is one of the most interesting questions
about the existence of the numerical Shilov boundary. In addition,
we show the existence of the numerical Shilov boundary for a locally
uniformly convex separable Banach space $X$.

\section{Differentiability of a subspace of $C_b(K)$}

\begin{defn}\label{defn:defnapproach}Let $K$ be a Hausdorff space and $A$ be a subspace of $C_b(K)$.
We say that {\it every norming sequence of $f$ approaches for $A$},
whenever for any two sequences $\{x_n\}_{n=1}^\infty$ and
$\{y_n\}_{n=1}^\infty$ in $K$ satisfying
\begin{equation}\label{eq:aapproach1}
\lim_n \alpha f(x_n) = \|f\|= \lim_n \beta f(y_n)
\end{equation}
for some complex numbers $\alpha$, $\beta$ of modulus $1$, we have
$\lim_n (\alpha g(x_n) - \beta g(y_n))=0$ for every $g\in A$. In
case that $\lim_n (\alpha g(x_n) - \beta g(y_n))=0$ uniformly for
$g\in S_A$, we say that {\it every norming sequence of $f$
approaches uniformly for $A$}.
\end{defn}

It is easy to see that if every norming sequence of $f$ approaches
for $A$, and also if $A$ is nontrivial, then $f\neq 0$.

\begin{thm}\label{thm:criterionGateaux}
Let $A$ be a subspace of $C_b(K)$and $f\in A$.
\begin{itemize}
\item[(i)] The norm $\|\cdot \|$ of $A$ is G\^{a}teaux differentiable
at $f$ if and only if every norming sequence of $f$ approaches for
$A$.

\item[(ii)] The norm $\|\cdot \|$ of $A$ is Fr\'echet differentiable at
$f$ if and only if every norming sequence of $f$ approaches
uniformly for $A$.
\end{itemize}
\end{thm}
\begin{proof} A slight modification of the proof of (i) gives that
of (ii), so we prove only (i). We may assume that $A$ is nontrivial.
Assume that $\|\cdot\|: A\to \mathbb{R}$ is G\^{a}teaux
differentiable at $f\in S_A$.  Take two sequences
$\{x_n\}_{n=1}^\infty$ and $\{y_n\}_{n=1}^\infty$ in $K$ satisfying
\[ \lim_n \alpha f(x_n) =\lim_n \beta f(y_n)=1=\|f\| \] for some
complex numbers $\alpha$, $\beta$  of modulus one. Fix $g\in A$.
Since the norm of $A$ is G\^{a}teaux differentiable at $f$, for
every $\epsilon>0$ there is $\delta>0$ such that \[ \|f + tg\| +
\|f-t g\|\le 2 + \epsilon |t|,\] for every real number $t$,
$|t|<\delta$. For every positive integer $n$ we have
\[ |f(x_n) + tg(x_n) | + |f(y_n)-t g(y_n)| \le 2 + \epsilon |t|,\]
and so
\[ \RE (\alpha f(x_n) + t \alpha g(x_n)) + \RE ( \beta f(y_n) - t
\beta g(y_n)) \le 2 + \epsilon |t|.\] Therefore, if $|t|<\delta$,
\[ \limsup_n t \RE (\alpha g(x_n) - \beta g(y_n)) \le \epsilon|t|,\]
 which implies that $\lim_n \RE(\alpha g(x_n) - \beta g(y_n))=0$.
 Replacing $g$ by $-ig$, we get $\lim_n \IM(\alpha g(x_n) -\beta
g(y_n))=0.$ Therefore, $\lim_n (\alpha g(x_n) -\beta g(y_n))=0$ for
every $g\in A.$

For the converse, assume that there is an $f\in S_A$ such that every
norming sequence of $f$ approaches for $A$, but $\|\cdot\|$ is not
G\^{a}teaux differentiable at $f$. Then there exist $g\in S_A$, a
null sequence $\{t_n\}$ of nonzero real numbers and $\epsilon>0$
such that
\begin{equation}\label{eq:nearnorm10} \|f+ t_n g\| + \|f-t_n g\|\ge
2+\epsilon |t_n|,\ \ \ \ \ \ \ \ \forall n\ge 1.\end{equation}

Choose sequences $\{x_n\}_{n=1}^\infty$ and $\{y_n\}_{n=1}^\infty$
in $K$ such that for each $n\ge 1$,
\begin{equation}\label{eq:nearnorm20} |(f+t_n g)(x_n)| \ge \|f+t_ng\|-\frac 1n |t_n|, \ \
|(f-t_ng)(y_n)| \ge \|f-t_ng\|-\frac 1n |t_n|.\end{equation} Then
\[ 1\ge |f(x_n)| \ge |(f+ t_ng)(x_n)| - |t_ng(x_n)|\ge \|f+t_ng\| -
 \frac 1n|t_n| - |t_ng(x_n)|.\] So it is clear that $\lim_n |f(x_n)|=1$.
Similarly, $\lim_n |f(y_n)|=1$.

Since every norming sequence of $f$ approaches for $A$, by passing
to a proper subsequence, we may assume that there exist two
sequences $\{x_n\}$, $\{y_n\}$ and complex numbers $\alpha$, $\beta$
of modulus one such that
\begin{equation}\label{eq:aapproesti210}
\lim_n \alpha f(x_n) = \lim_n \beta f(y_n) = 1\text{ and }\sup_{n\ge
1}|\alpha g(x_n) - \beta g(y_n)| \le \epsilon/2.
\end{equation}

Using \eqref{eq:nearnorm10}, \eqref{eq:nearnorm20} and
\eqref{eq:aapproesti210}, we get for any $n$,

\begin{align*}
2 + \epsilon |t_n| &- \frac{2}n |t_n|\le  \|f+t_ng\| +
\|f-t_ng\|-\frac{2}n |t_n|\\&\le |(f+t_ng)(x_n)| + |(f-t_ng)(y_n)|
\\&=|(\alpha f+\alpha t_ng)(x_n)| + |(\beta f-\beta t_ng)(y_n)|\\
&\le |(\alpha f+\alpha t_ng)(x_n)| + |\beta f(y_n)-\alpha t_ng(x_n)|
+ |\alpha t_ng(x_n) -\beta t_ng(y_n)|
\\ &\le |\alpha f(x_n)+\alpha t_ng(x_n)| + |\beta f(y_n)-\alpha t_ng(x_n)| +
\frac \epsilon2 |t_n|
\end{align*}
Hence for every $n\ge 1$,
\begin{equation}\label{eq:basicestimate110}
2 + (\frac \epsilon2 -\frac{2}{n})|t_n| \le |\alpha f(x_n)+\alpha
t_ng(x_n)| + |\beta f(y_n)-\alpha t_ng(x_n)|.
\end{equation}

We need the following basic lemma which is proved later.
\begin{lem}\label{lem:basiccalculus}
Let $\varphi:U\subset \mathbb{R}^n\to \mathbb{R}$ be twice
continuously differentiable on a neighborhood $U$ of $\xi_0\in
\mathbb{R}^n$. Let $\epsilon>0$. Then there exist $\delta>0$ and a
neighborhood $V$ of $\xi_0$ such that for any $\xi, \zeta$ in $V$
and $|\eta|<\delta$,
\[
|\varphi(\xi+\eta) - \varphi(\xi) + \varphi(\zeta-\eta) -
\varphi(\zeta)| \le \epsilon|\eta|.
\]
\end{lem}

Notice that the function $\varphi:\mathbb{R}^2 \to \mathbb{R}$
defined by $\varphi(\xi) = |\xi|$ is infinitely differentiable on a
neighborhood of $(1,0)$, where $|\cdot|$ is a usual Euclidean norm
in $\mathbb{R}^2$. By Lemma~\ref{lem:basiccalculus}, given
$\epsilon>0$, there exist a neighborhood $V$ of $(1,0)$ and a
$\delta>0$ such that for any $\xi, \zeta$ in  $V$ and
$|\eta|<\delta$,
\begin{equation}\label{eq:basicestimate210} |\varphi(\xi+\eta) -
\varphi(\xi) + \varphi(\zeta-\eta) - \varphi(\zeta)| \le
\epsilon|\eta|/4.\end{equation}

We shall identify the complex plane $\mathbb{C}$ with
$\mathbb{R}^2$. For each $n$, set $\xi_n = \alpha f(x_n)$, $\zeta_n
= \beta f(y_n)$ and $\eta_n = \alpha t_ng(x_n)$. By
\eqref{eq:aapproesti210}, we may assume that $\xi_n$ and $\zeta_n$
are in $V$ and $|\eta_n|<\delta$ for any $n$. By
\eqref{eq:basicestimate210}, for every $n$,
\begin{align}\label{eq:basicestimate310}
|\alpha f(x_n) + \alpha t_ng(x_n)|- |f(x_n)| &+ |\beta f(y_n) -
\alpha t_ng(x_n)| - |f(y_n)|\nonumber
\\&\le |\varphi(\xi_n+\eta_n) - \varphi(\xi_n) + \varphi(\zeta_n-\eta_n) -
\varphi(\zeta_n)|
\\&\le  \epsilon|\eta_n|/4 =\epsilon|t_ng(x_n)|/4\le
\epsilon|t_n|/4.\nonumber
\end{align}

By \eqref{eq:basicestimate110} and \eqref{eq:basicestimate310}, we
get for every $n$,
\[   2 + (\frac \epsilon2
-\frac{2}{n})|t_n|- |f(x_n)|-|f(y_n)| \le \frac{\epsilon}4|t_n|.
\]
This means that\[ 0\le 2-(|f(x_n)|+|f(y_n)|) \le (-\frac{\epsilon}4
+ \frac 2n) |t_n|<0\] for sufficiently large $n$. This is a
contradiction. The proof is done.
\end{proof}

Now we prove Lemma~\ref{lem:basiccalculus}.
\begin{proof}[Proof of Lemma~\ref{lem:basiccalculus}]
Choose a positive $r>0$ such that $B(\xi_0, 4r)=\{ \xi: |\xi-\xi_0|
\le 4r\}$ is contained in $U$. For any $\xi \in B(\xi_0, r)$ and
$|\eta|<r$, by the Taylor formula of $\varphi$ there is $0\le t\le
1$ such that
\[
\varphi(\xi + \eta) - \varphi(\xi)- \nabla \varphi(\xi)\cdot \eta
=\frac 12\sum_{i,j=1}^n \frac{\partial^2 \varphi}{\partial
\xi_i\partial \xi_j}(\xi+t\eta)\eta_i \eta_j,
\] where $\nabla \varphi(\xi)\cdot \eta =\sum_{i=1}^n \frac{\partial
\varphi}{\partial \xi_i}(\xi)\eta_i$. Let $M=\sup_{\xi\in B(\xi_0,
2r)} \frac{\partial^2 \varphi}{\partial \xi_i\partial \xi_j}(\xi).$
Then for $\xi\in B(\xi_0, r)$ and $|\eta|<r$,
\begin{equation}\label{eq:difference11000}
|\varphi(\xi + \eta) - \varphi(\xi) - (\nabla\varphi)(\xi)\cdot \eta
|\le \frac12 n^2M|\eta|^2\end{equation}

Notice that the mapping $\xi \to \nabla\varphi(\xi)$ from $B(\xi_0,
4r)$ to $\mathbb{R}^n$ is uniformly continuous. By
\eqref{eq:difference11000}, given $\epsilon >0$ there exists
$\delta>0$ such that for any $\xi, \zeta$ in $B(\xi_0, \delta)$ and
for any $|\eta|<\delta$,
\[
|\varphi(\xi + \eta) - \varphi(\xi) - (\nabla\varphi)(\xi)\cdot \eta
|\le \epsilon |\eta|/4
\] and
\[ |\nabla \varphi(\xi) - \nabla \varphi(\zeta)|\le\epsilon/2.\]

Take $V=B(\xi_0, \delta)$. For any $\xi, \zeta$ in $V$ and
$|\eta|<\delta$,
\begin{align*}
|\varphi(\xi+\eta) - &\varphi(\xi) + \varphi(\zeta-\eta) -
\varphi(\zeta)| \\ &\le |\varphi(\xi+\eta) - \varphi(\xi) -
\nabla\varphi(\xi)\cdot \eta | + |\varphi(\zeta-\eta) -
\varphi(\zeta) + \nabla\varphi(\zeta)\cdot \eta|
\\&\ \ \ \ \ \ +|\nabla\varphi(\xi)\cdot \eta - \nabla\varphi(\zeta)\cdot \eta|
\\& \le \epsilon|\eta|/2 +
|\nabla\varphi(\xi)-\nabla\varphi(\zeta)|\cdot|\eta|\le
\epsilon|\eta|.
\end{align*}
The proof is done.\end{proof}

Notice that $X$ and $X^*$ can be regarded as a subspace of
$C(B_{X^*})$ and $C_b(B_{X})$ respectively, where the weak-$*$ and
norm topology is given on $B_{X^*}$ and  $B_{X}$, respectively. By
the direct application of Theorem~\ref{thm:criterionGateaux} we get
the following \v{S}mulyan's theorem.

\begin{thm}[\v{S}mulyan]\label{thm:Smulyan}
Let $X$  be a Banach space. Then
\begin{itemize}
\item[(i)] The norm of $X$ is Fr\'echet differentiable at $x\in
S_X$ if and only if whenever $x_n^*, y_n^*\in S_{X^*}$, $x_n^*(x)\to
1$ and $y_n^*(x)\to 1$, then $\|x_n^* - y^*_n\|\to 0$.

\item[(ii)] Then norm of $X^*$ is Fr\'echet differentiable at $x^*\in
S_{X^*}$ if and only if whenever $x_n, y_n\in S_X$, $x^*(x_n)\to 1$
and $x^*(y_n)\to 1$, then $\|x_n - y_n\|\to 0$.

\item[(iii)] The norm of $X$ is G\^ateaux differentiable at $x\in
S_X$ if and only if whenever $x^*_n, y^*_n\in S_{X^*}$, $x^*_n(x)\to
1$ and $y^*_n(x)\to 1$, then
$x_n^*-y^*_n\stackrel{w^*}{\longrightarrow} 0.$

\item[(iv)] The norm of $X^*$ is G\^ateaux differentiable at $x^*\in
S_{X^*}$ if and only if $x_n, y_n\in S_X$, $x^*(x_n)\to 1$ and
$x^*(y_n)\to 1$, then $x_n - y_n \stackrel{w}{\longrightarrow} 0.$
\end{itemize}
\end{thm}

\begin{prop}\label{prop:strongmeansgateaux}
Let $K$ be a Hausdorff space and $A$ be a subspace of $C_b(K)$.
\begin{itemize}
\item[(i)] If $f$ is a strong peak function in $A$, then every norming
sequence of $f$ approaches for $A$. Hence the norm $\|\cdot \|$ is
G\^{a}teaux differentiable at every strong peak function.

\item[(ii)] Assume in addition that $A$ is a separating subspace of
$C(K)$ on a compact Hausdorff space $K$ and that $f$ is a nonzero
element of $A$. Then the norm of $A$ is G\^{a}teaux differentiable
at $f$ if and only if $f$ is a strong peak function. In this case,
the set of all weak-$*$ exposed points of $B_{A^*}$ is
\[ w^{\!*\!}exp B_{A^*} = \{ \alpha \delta_t\ :\ t\in \rho A,\  |\alpha|=1\}.\]
\end{itemize}
\end{prop}
\begin{proof}
(i) Suppose that $f$ is a strong peak function at $x_0$ and that
there exist two sequences $\{x_n\}_{n=1}^\infty$ and
$\{y_n\}_{n=1}^\infty$ in $K$ satisfying
\[ \lim_n \alpha f(x_n) = \lim_n \beta f(y_n) = \|f\|\] for some
complex numbers $\alpha$, $\beta$ of modulus one. Then, two
sequences converge to $x_0$ in $K$ and $\alpha=\beta$. It is clear
that $\lim_n (\alpha g(x_n) -\beta g(y_n))=0$ for every $g\in A$.
This completes the proof of (i).

(ii) It is enough to prove the necessity. We may assume $\|f\|=1.$
Since the norm $\|\cdot \|$ is G\^{a}teaux differentiable at $f$,
$f$ is a smooth point of $B_A$. Choose $t\in K$ and $\alpha,
~|\alpha|=1$ such that $\alpha f(t) = 1$. Then the evaluation
functional $\alpha \delta_t\in S_{A^*}$ is a weak-$*$ exposed point
of $B_{A^*}$. Since $A$ is separating, $\alpha \delta_t \neq \beta
\delta_s$ on $A$ if $t\neq s$ in $K$ and $\alpha, \beta\in
 S_\mathbb{C}$. If $s\neq t$,
\[ \|f\|=|f(t)|=1> \max \{ \RE \beta f(s): \beta\in S_\mathbb{C}\} =
|f(s)|.\] Hence $f$ is a peak function in $A$.
\end{proof}

Consider the product space $K\times B_{Y^*}$, where $B_{Y^*}$ is
equipped with the weak-$*$ topology. Given a subspace $A$ of
$C_b(K,Y)$, consider the map $\varphi: f\in A \mapsto \tilde{f}\in
C_b(K\times B_{Y^*})$ defined by
\[ \tilde{f}(x, y^*) = y^*f(x), \ \ \forall (x, y^*)\in K\times B_{Y^*}.\]
Then $\varphi$ is a linear isometry, and its image $\tilde{A}$ of
$A$ is also a subspace of $C_b(K\times B_{Y^*})$. In particular, we
shall say that the subspace $A$ of $C_b(K, Y)$ is {\it separating}
if the following conditions hold:

\begin{itemize}
\item[(i)] If $x\neq y$ in $K$, then $\delta_{(x,x^*)} \neq
\delta_{(y,y^*)}$ on $\tilde{A}$ for $x^*, y^* \in S_{Y^*}$.

\item[(ii)] Given $x\in K$ with $\delta_x \neq 0$ on $A$, we have
$\delta_{(x,x^*)} \neq \delta_{(x,y^*)}$ on $\tilde{A}$ for $x^*\neq
y^*$ in $ext(B_{Y^*})$.
\end{itemize}
By applying Theorem~\ref{thm:criterionGateaux} and
Proposition~\ref{prop:strongmeansgateaux} to the subspace
$\tilde{A}$ of $C_b(K\times B_{Y^*})$, we get the following.

\begin{cor}\label{cor:vectorGateaux}
Let $K$ be a Hausdorff space, $Y$ a Banach space and $A$ a subspace
of $C_b(K, Y)$. Then the following hold:
\begin{itemize}
\item[(i)] \label{cor:vectorGateaux1}The norm $\|\cdot
\|$ of $A$ is G\^{a}teaux differentiable (resp. Fr\'echet
differentiable) at $f$ if and only if whenever there exist sequences
$\{x_n\}_{n=1}^\infty$, $\{y_n\}_{n=1}^\infty$ in $K$ and
$\{x^*_n\}_{n=1}^\infty$, $\{y_n^*\}_{n=1}^\infty$ in $S_{Y^*}$ such
that
\[ \lim_{n\to \infty} x^*_nf(x_n) = \|f\| = \lim_{n\to \infty}
y_n^*f(y_n),\] we get
\[  \lim_{n\to \infty} (x^*_ng(x_n) -y_n^*g(y_n))=0,\ \ \ \ \forall g\in A.\]
\[ (\text{ resp. } \lim_{n\to \infty} (x^*_ng(x_n) -y_n^*g(y_n))=0\ \text{uniformly for }\   g \in S_A \ ).\]

\item[(ii)] \label{cor:vectorGateaux3}If $f$ is a strong peak function at $x_0\in K$ and $f(x_0)/\|f(x_0)\|_Y$ is a smooth point of $B_Y$,
then the norm $\|\cdot \|$ of $A$ is G\^{a}teaux differentiable at
$f$.

\item[(iii)] \label{cor:vectorGateaux4}Assume in addition that $A$ is a separating subspace of
$C_b(K, Y)$ on a compact Hausdorff space $K$ and that $f$ is a
nonzero element of $A$. Then the norm of $A$ is G\^{a}teaux
differentiable at $f$ if and only if $f$ is a strong peak function
at some $x_0$ and $f(x_0)/\|f(x_0)\|_Y$ is a smooth point of $B_Y$.
In this case, the set of all weak-$*$ exposed points of $B_{A^*}$ is
\begin{align*}
w^{\!*\!}exp B_{A^*} = \{ \alpha \delta_{(x, y^*)}\ :~ \exists &
\text{ a strong peak function }  f \text{ such that } \\ &y^*f(x)
=\|f\|  \text{ and } \ y^*\in w^*exp(B_{Y^*})\ \},
\end{align*} where $\delta_{(x, y^*)}(f)= y^*f(x)$ for all $f\in
C_b(K, Y)$.

\end{itemize}
\end{cor}

\begin{prop}
If we denote by $G'(f)$ the G\^{a}teaux differential of the norm at
$f$, then \[ G'(f)(g) = \lim_n  \RE(\alpha g(x_n)),\] for a sequence
$\{x_n\}_n$ in $K$ and a complex number $\alpha$ of modulus $1$
satisfying $\lim_n \alpha f(x_n) = \|f\|$.
\end{prop}
\begin{proof}
Since $\RE(\alpha f(x_n) + \alpha tg(x_n))\le \|f+tg\|$, we have
$t\RE(\alpha g(x_n)) \le \|f+ tg\| - \RE(\alpha f(x_n))$ for all
real $t$. Hence for $t>0$,
\[ \limsup_n \RE(\alpha g(x_n)) \le \lim_n \frac{\|f+tg\|-\RE(\alpha
f(x_n))}{t}=\frac{\|f+tg\|-\|f\|}{t},\] and  for $t<0$,
\[ \lim\inf_n \RE(\alpha g(x_n)) \ge \lim_n \frac{\|f+tg\|-\RE(\alpha
f(x_n))}{t}=\frac{\|f+tg\|-\|f\|}{t}.\] Therefore, it is easy to see
that
\[ \lim_n \RE(\alpha g(x_n)) = \lim_{t\to0}
\frac{\|f+tg\|-\|f\|}{t}= G'(f)(g).\] This completes the proof.
\end{proof}

We apply Theorem~\ref{thm:criterionGateaux} to show that the norm of
$A(B_X)$ is nowhere Fr\'echet differentiable, if $X$ is nontrivial.

\begin{prop}\label{prop:strongpeaknotuniform}
Suppose that $X$ is a nontrivial complex Banach space and that $f$
is a strong peak function in $A(B_X)$. Then every norming sequence
of $f$ doesn't approach uniformly for $A(B_X)$.
\end{prop}
\begin{proof}
Let $f\in A(B_X)$ be a strong peak function at some $x_0\in S_X$
After a proper rotation, we may assume that $f(x_0)=\|f\|$. Let $x_n
= e^{i/n} x_0$ for every positive integer $n$. It is easy to see
that $\{x_n\}$ is a norming sequence and each $x_n$ is a strong peak
point for $A(B_X)$, so there is a strong peak function $g_n\in
A(B_X)$ such that $g(x_n)=1= \|g_n\|$ and $|g_n(x_0)|<1/2$. Hence we
get for every $n$,
\begin{equation}\label{eq1:nowherediff}
|g_n(x_n)-g_n(x_0)|\ge |g_n(x_n)| - |g_n(x_0)|\ge 1/2.\end{equation}
Since $\lim_n f(x_n)=\|f\| = f(x_0)$, \eqref{eq1:nowherediff}
implies that every norming sequence of $f$ doesn't approach
uniformly for $A(B_X)$.
\end{proof}

\begin{thm}Suppose that $X$ is a nontrivial complex Banach
space. The norm $\|\cdot \|$ of $A(B_X)$ is nowhere Fr\'echet
differentiable.
\end{thm}
\begin{proof}
Suppose that the norm of $A(B_X)$ is Fr\'echet differentiable at
some $f$. By Theorem~ \ref{thm:criterionGateaux}, every norming
sequence of $f$ approaches uniformly for $A(B_X)$. This implies that
$f$ is a strong peak function which contradicts
Proposition~\ref{prop:strongpeaknotuniform}. In fact, suppose that
there is a sequence $\{x_n\}$ in $S_X$ such that $$\lim
|f(x_n)|=\|f\|.$$ By passing to a proper subsequence, we may assume
that there is a complex number $\alpha,~|\alpha|=1$ such that $\lim
\alpha f(x_n) =\|f\|$. We claim that the sequence $\{x_n\}$ is
Cauchy. Otherwise, there exist subsequence $\{x_{n_k}\}$  and
$\delta>0$ such that for $\|x_{n_{k+1}} - x_{n_k}\|\ge \delta$ for
every $k$. Since \[\lim \alpha f(x_{n_{k}}) = \|f\| = \lim \alpha
f(x_{n_{k+1}}),\] and since every norming sequence of $f$ approaches
uniformly for $A(B_X)$, we have $\lim_k |g(x_{n_{k+1}})-
g(x_{n_k})|=0$ uniformly in $g\in S_{A(B_X)}$. Since $S_{X^*}\subset
S_{A(B_X)}$, we have that $\lim_k \|x_{n_{k+1}}- x_{n_k}\|=0$, which
is a contradiction. Let $x_0$ be a limit of $\{x_n\}$. Suppose that
there is another sequence $\{y_n\}$ in $S_X$ such that $\lim_n
|f(y_n)|=\|f\|$. By choosing an appropriate subsequence, we may
assume that there is a complex number $\beta,~|\beta|=1$ such that $
\lim_n \beta f(y_n) = \|f\|$. Then $ \lim_n \beta f(y_n) = \|f\| =
\lim_n \alpha f(x_n)$. Since every norming sequence of $f$
approaches uniformly for $A(B_X)$, $\alpha=\beta$ and  $\lim_n \|
x_n - y_n\|=0$. Therefore, $\lim_n y_n = \lim_n x_n =x_0$. This
shows that $f$ is a strong peak function at $x_0$.\end{proof}

\begin{remark}
When $X=\{0\}$, it is easy to see that $A(B_X)$ is isometrically
isomorphic to $\mathbb{C}$. Thus the norm is Fr\'echet
differentiable everywhere except zero.
\end{remark}

\section{Bishop's theorem}

Bishop  showed in \cite{B} that if $K$ is a compact metrizable and
if $\rho A$ is the set of all (strong) peak points for a separating
function algebra $A$, then
\[ \max_{t\in \rho A} |f(t)| = \|f\| \ \ \ \ \text{for every } f \in
A.\] We now give another version of Bishop's theorem from the
results in the previous section.

\begin{thm}[Bishop's theorem]\label{thm:Bishop}
Let $A$ be a nontrivial separating separable subspace of $C(K)$ on a
compact Hausdorff space $K$. Then the set of all peak functions in
$A$ is a dense $G_\delta$-subset of $A$. In particular, $\rho A$ is
a norming subset for $A$ and  $\partial A=\overline{\rho A}$.
\end{thm}
\begin{proof}
By Proposition~\ref{prop:strongmeansgateaux} and Mazur's theorem,
the set of all peak functions in $A$ is a dense $G_\delta$-subset of
$A$. It is clear that every closed boundary for $A$ contains $\rho
A$. Hence we have only to show that $\rho A$ is a norming subset for
$A$. For each $f\in A$, there is a sequence $\{f_n\}$ of peak
functions such that $\|f_n - f\|\to 0$ as $n\to \infty$. Then
\begin{align*}-\|f-f_n\| + \|f_n\|&\le -|f(x_n) - f_n(x_n)| +
|f_n(x_n)|\le |f(x_n)| \\&\le |f(x_n) - f_n(x_n)| + |f_n(x_n)|\le
\|f-f_n\| + \|f_n\|,\end{align*} where $x_n$ is a peak point for
$f_n$ for each $n$. Hence $\lim_n |f(x_n)| = \lim_n \|f_n\| =
\|f\|.$ Notice that $x_n\in \rho A$. Therefore, $\rho A$ is a
norming subset for $A$. The proof is done.
\end{proof}

The following example given in \cite{B} shows that the separability
assumption in Theorem~\ref{thm:Bishop} is necessary. Let $J$ be an
uncountable set and let $I_\alpha = [0,1]$ for each $\alpha\in J$.
Then the product space $K= \Pi_{\alpha} I_\alpha$ is a compact
non-metrizable space, and  $C(K)$ is not separable. Using the
Stone-Weierstrass theorem, it is not difficult to check that for
every $f\in C(K)$, there is a countable subset $\Delta\subset J$
such that whenever $x$ and $y$ in $K$ satisfy $x_\alpha = y_\alpha$
for every $\alpha \in \Delta$,  $f(x) = f(y)$ holds. It is easy to
see that there is no peak function in  $C(K)$. In particular, the
norm of $C(K)$ is nowhere G\^{a}teaux differentiable by
Proposition~\ref{prop:strongmeansgateaux}.

\begin{example}
Let $X=\ell_2^2$ be the 2-dimensional complex Euclidean space, and
let $A$ be the set of restrictions to $B_X$ of the elements of
$X^*$, which is a closed subspace of $C(B_X)$. Given two distinct
points $x, y\in B_X$, there is $f\in A$ such that $f(x)\neq f(y)$,
but it is easy to see that the subspace $A$ is not separating. The
set $T_1=\{(x_1, x_2): (x_1, x_2)\in S_X,~ x_2\ge 0 \}$ and
$T_2=\{(x_1, -x_2): (x_1, x_2)\in S_X,~ x_2\ge 0 \}$ are two closed
norming subsets for $A$. However, $T_1\cap T_2= \{ (x_1, 0):
|x_1|=1\}$ is not a norming subset for $A$, so $A$ doesn't have the
Shilov boundary. Therefore we cannot omit the separation assumption
in Theorem~\ref{thm:Bishop}.
\end{example}

The following  is a consequence of Corollary~\ref{cor:vectorGateaux}
and Theorem~\ref{thm:Bishop}.
\begin{cor}\label{cor:vectorBishop}
Let $Y$ be a Banach space and let $A$ be a nontrivial separating
separable subspace of $C(K, Y)$ on a compact Hausdorff space $K$.
Then the set
\[\{ f\in A:f \text{ is a peak function at some } t\in K,
f(t)/\|f\| \text{ is a smooth point of } B_Y\}\] is a dense
$G_\delta$-subset of $A$. In particular, $\rho A$ is a norming
subset for $A$ and $\partial A =\overline{\rho A}$.
\end{cor}

Notice that if $K$ is a compact metric space and $Y$ is separable,
then every subspace $A$ of $C_b(K, Y)$ is separable. Indeed, we can
regard $A$ as a subspace of $C(K\times B_{Y^*})$ and $K\times
B_{Y^*}$ is a compact metrizable space.

From the proof of Theorem~\ref{thm:Bishop} we have the following
proposition.

\begin{prop}\label{prop:strongpeakfunctiondense}
Given a Banach space $Y$, let $A$ be a  nontrivial subspace of
$C_b(K,Y)$ on a Hausdorff space $K$. Suppose that the set of all
strong peak functions in $A$ is dense. Then $\rho A$ is a norming
subset for $A$ and  $\partial A =\overline{\rho A}$.
\end{prop}

Assume that $A$ is a subspace of $C(K)$ on compact Hausdorff space
$K$ and for any two distinct points $s,t$ in $K$, there is $f\in A$
such that $f(s)\neq f(t)$. Then the mapping $x\mapsto \delta_x$ from
$K$ into the weak-$*$ compact subset $B_{A^*}$ is an injective
homeomorphism and we shall identify $K$ with its image in $B_{A^*}$.

\begin{prop}Suppose that $A$ is a subspace of $C(K)$ on a compact
Hausdorff space $K$ and that for two distinct points $t,s\in K$,
there is $f\in A$ such that $f(t)\neq f(s)$. Then $A$ is separable
if and only if $K$ is metrizable.
\end{prop}
\begin{proof}
Recall that if $A$ is separable, then the weak-$*$ compact set
$B_{A^*}$ is metrizable. Since $K$ is embedded in $B_{A^*}$, it is
metrizable. For the converse, notice that if $K$ is metrizable, the
Stone-Weierstrass theorem shows that $C(K)$ is separable, and so is
its subspace $A$.
\end{proof}

\section{Density of strong peak functions in $A(B_X, Y)$.}

Let $C$ be a closed convex and bounded set in a Banach space $X$.
The set $C$ is said to have the {\it Radon-Nikod\'ym property} if
for every probability space $(\Omega, \mathcal{B}, \mu)$ and every
$X$-valued countably additive measure $\tau$ on $\mathcal{B}$ such
that $\tau(A)/\mu(A)\in C$ for every $A\in \mathcal{B}$ with
$\mu(A)>0$, there is a Bochner measurable $f:\Omega\to X$ such that
\[ \tau(A) = \int_A f(\omega) \, d\mu(\omega), \ \ \ \ A\in
\mathcal{B}.\]

A Banach space $X$ is said to have the {\it Radon-Nikod\'ym
property} if its unit ball $B_X$ has the Radon-Nikod\'ym property.
For the basic properties and useful information on the
Radon-Nikod\'ym property, see \cite{DU, FLP, JL}.

Let $D$ be a metric space. We say that a function $\varphi:D\to
\mathbb{R}$ {\it strongly exposes} $D$ if there is $x\in D$ such
that
\[ \varphi(x) = \sup\{\varphi(y) : y\in D\}\] and whenever there is
a sequence $\{x_n\}$ in $D$ satisfying $\lim_n \varphi(x_n) =
\varphi(x)$, the sequence $\{x_n\}$ converges to $x$.

The important Bourgain-Stegall's perturbed optimization theorem
\cite{S} says that if a closed bounded convex subset $D$ of $X$ has
the Radon-Nikod\'ym property and if $\varphi:D\to \mathbb{R}$ is a
bounded above upper semi-continuous function, then the set
\[ \{ x^* : \varphi + x^* \ \text{strongly exposes } D\}\] is a
dense $G_\delta$-subset of $X^*$.

Let $X$ and $Y$ be complex Banach spaces. Notice that $f\in
A(B_X,Y)$ is a strong peak function if and only if
 $\|f(\cdot)\|$ strongly exposes $B_X$.

\begin{defn}
A function $f\in A(B_X, Y)$ is said to {\it attain its norm strongly
} on $B_X$ if there is $x_0 \in S_X$ such that whenever $\lim_n
\|f(x_n)\|=\|f\|$ for a sequence $\{x_n\}$ in $B_X$, it has a
subsequence $\{x_{n_k}\} $ converging to $\alpha x_0$ for some
$|\alpha|=1$.
\end{defn}

Acosta, Alaminos, Garc\'{i}a and Maestre \cite{AAGM} showed that if
$X$ has the Radon-Nikod\'{y}m property, then for every $f\in A(B_X,
Y)$, every natural number $N$ and every $\epsilon>0$, there are
$x_1^*,\ldots, x_N^*\in X^*$ and $y_0\in Y$ such that the
$N$-homogeneous polynomial $Q$ on $X$, given by $Q(x)=
x_1^*(x)\cdots x_N^*(x) y_0$ satisfies that $\|Q\| < \epsilon$ and
$f+Q$ attains its norm. For our application we prove the following
stronger version.

\begin{thm}\label{thm:radonstrongnormattaining}
Let $X$ be a complex Banach space with the Radon-Nikod\'ym property.
Suppose that $f\in A(B_X, Y)$, $N\ge 1$ and $\epsilon>0$. Then there
are $x_1^*,x_2^*\in X^*$ and $y_0\in Y$ such that the
$N$-homogeneous polynomial $Q$ on $X$, given by $Q(x)=
[x_1^*(x)]^{N-1} x_2^*(x) y_0$, satisfies that $\|Q\| < \epsilon$
and $f+Q$ strongly attains its norm. In particular, the set of all
strongly norm-attaining functions is dense in $A(B_X, Y)$.

\end{thm}
\begin{proof}
We may assume that $X\neq 0$. Fix $f\in A(B_X, Y)$ and define a
function $g:B_X\to \mathbb{C}$ as the following:
\begin{equation}\label{eq:symmetricform} g(x) = \max\{ \|f(\lambda
x)\| : \lambda\in \mathbb{C}, |\lambda|\le 1\}.\end{equation} It is
clearly bounded, because $f$ is an element of $A(B_X,Y)$.

For the proof of the upper semi-continuity of $g$, suppose that a
sequence $\{x_n\}_{n=1}^\infty$ in $B_X$ converges to $x$. Then for
each $n$, there is a complex number $\lambda_n$ such that
$|\lambda_n|=1$ and $g(x_n) =\|f(\lambda_n x_n)\|$. For any
convergent subsequence $\{\lambda_{n_k}\}$ of $\{\lambda_n\}$ with
the limit $\lambda$, we get
\[\lim_{k\to \infty}\|f(\lambda_{n_k} x_{n_k})\|= \|f(\lambda x)\|\le
g(x).\] Hence $\limsup_n g(x_n)\le g(x)$. This means that $g$ is
upper semi-continuous.

By Bourgain-Stegall's perturbed optimization theorem, there is
$x^*\in X^*$ such that $\|x^*\|<\epsilon$ and $g + \RE x^*$ strongly
exposes $B_X$ at $x_0$.

We claim that $\RE x^*(x_0)\neq 0$. Assume that $\RE x^*(x_0)=0$.
Then $g(x_0)+\RE x^*(x_0)= g(-x_0) + \RE x^*(-x_0)$. So $x_0=0$.
Notice that for each $x\in B_X$, $g(0)=\|f(0)\|\le \|f(x)\|\le g(x)$
by the maximum modulus theorem. Since $g+\RE x^*$ strongly exposes
$B_X$ at 0,
\begin{align*} g(0) &= \sup\{ g(x) + \RE x^*(x) : x\in B_X\}\\
&= \sup\{ g(x) +  |x^*(x)| : x\in B_X\}.
\end{align*} Hence $g(0)\le g(x)\le g(x)+ |x^*(x)|\le g(0)$ for any
$x\in B_X$. This means that $x^*=0$ and $g$ is constant on $B_X$.
This is a contradiction to that $g$ strongly exposes $B_X$ at 0.
Therefore $\RE x^*(x_0)\neq 0$.

Then $\|x_0\|=1$. Indeed, it is clear that $x_0\neq 0$, because
$x^*\neq0$ and $g$ is nonnegative. If $0<\|x_0\|<1$, then
\begin{align*} g(x_0) + \RE
x^*(x_0) &= \sup\{ g(x) + \RE x^*(x) : x\in B_X\}\\
&= \sup\{ g(x) +  |x^*(x)| : x\in B_X\}
\end{align*}
shows that $\RE x^*(x_0) = |x^*(x_0)|$ and
\[ g(x_0) + |x^*(x_0)| < g(\frac{x_0}{\|x_0\|}) +
|x^*(\frac{x_0}{\|x_0\|})| = g(\frac{x_0}{\|x_0\|}) + \RE
x^*(\frac{x_0}{\|x_0\|}).\] This is a contradiction to the fact that
$g +\RE x^*$ strongly exposes $B_X$ at $x_0$.

There is a $\lambda_0$ such that $|\lambda_0|=1$ and $g(x_0)=
\|f(\lambda_0 x_0)\|$. Let $x_1= \lambda_0 x_0$ and choose $x^*_1\in
X^*$ with $x_1^*(x_1)=1=\|x^*_1\|.$ Define $h:B_X\to Y$ by
\[ h(x) = f(x) + \lambda_1 x^*_1(x)^{N-1} x^*(x)\frac{f(x_1)}{\|f(x_1)\|},\] where
the complex number $\lambda_1$ is  properly chosen so that
\[ \big|\|f(x_1)\| +\lambda_1 x^*(x_1)\big| = \|f(x_1)\|+|x^*(x_1)|.\]

It is clear that $h\in A(B_X, Y)$ and notice that we get for every
$x\in B_X$, \begin{align} \|h(x)\| &\le \|f(x)\| + |x^*(x)| \le g(x)
+ |x^*(x)| \label{eq13}\\&\le \sup\{ g(x) + |x^*(x)| : x\in B_X\}
\notag\\ &\ = \sup\{ g(x) +  \RE x^*(x) : x\in B_X\}=g(x_0) + \RE
x^*(x_0).\notag\end{align} Hence $\|h\| = g(x_0) + \RE x^*(x_0)$
because $\RE x^*(x_0) = |x^*(x_0)|$ and
\begin{align*}\|h(x_1)\|& = |\|f(x_1)\|+\lambda_1 x^*(x_1)|= \|f(x_1)\|+
|x^*(x_1)|\\& = g(x_0) + |x^*(x_0)| = g(x_0) + \RE
x^*(x_0).\end{align*}

We shall show that $h$ strongly attains its norm at $x_0$. Suppose
that $\lim_n \|h(x_n)\|=\|h\| = g(x_0) + \RE x^*(x_0)$. Choose a
sequence  $\{\alpha_n\}$ of complex numbers so that $|\alpha_n|=1$
and
\[ g(x_n) + |x^*(x_n)| = g(\alpha_n x_n) + \RE x^*(\alpha_n x_n),\ \ \ \forall n\ge 1.\]
Then \eqref{eq13} shows that
\[\lim_{n\to \infty}g(\alpha_n x_n) + \RE x^*(\alpha_n x_n) = g(x_0) + \RE
x^*(x_0).\] Since $g+ \RE x^*$ strongly exposes $B_X$ at $x_0$,
$\{\alpha_n x_n\}$ converges to $x_0$. Hence there is a subsequence
of $\{x_n\}$ which converges to $\alpha x_0$ for some $|\alpha|=1$.
This implies that $h$ strongly attains its norm at $x_0$ and
$\|f-h\|\le \epsilon$. The proof is done.
\end{proof}

\begin{remark}In \eqref{eq:symmetricform} $g$ is continuous,
because it is the supremum of a family of continuous functions, that
is, it is lower semi-continuous.
\end{remark}

\begin{thm}\label{thm:densitypeakRNP}
Suppose that a complex Banach space $X$ has the Radon-Nikod\'ym
property and $Y$ is a nontrivial complex Banach space. Then the
following hold:
\begin{itemize}
\item[(i)] The set of all strong peak functions in $A(B_X, Y)$ is dense
in $A(B_X, Y)$. In particular, the set of all smooth points of
$B_{A(B_X, Y)}$ is dense in $S_{A(B_X, Y)}$ if the Banach space $Y$
is smooth.

\item[(ii)] $\rho A(B_X)$ is a norming subset for $A(B_X,
Y)$, and $\partial A(B_X,Y)=\overline{\rho A(B_X)}$. In particular,
$ext_\mathbb{C}(B_X)$ is a norming subset for $A(B_X, Y)$.
\end{itemize}
\end{thm}
\begin{proof} (i)  Suppose that $f\in A(B_X,Y)$ strongly attains its norm at $x_0$. We shall
show that given $\epsilon >0$ there is  $\tilde{g}\in A_u(B_X, Y)$
such that $\|\tilde{g}\|\le \epsilon$ and $f+\tilde{g}$ is a strong
peak function in $A(B_X, Y)$.

Since $f$ strongly attains its norm at $x_0$, there is a complex
number $\alpha$ of modulus $1$ such that $\|f(\alpha x_0)\| =
\|f\|$. Choose $x^*\in S_{X^*}$ so that $x^*(x_0)=1$ and take a peak
function $g\in A_u(\mathbb{D})$ such that $g(\alpha)=1$ and
$|g(\gamma)|<1$ for any $\gamma\neq \alpha$, where $\mathbb{D}$ is
the closed unit disc in $\mathbb{C}$. Define $h:B_X\to Y$ by
\begin{equation}\label{eq:peak1} h(x)= f(x)+\epsilon
g(x^*(x))\frac{f(\alpha x_0)}{\|f(\alpha x_0)\|}.\end{equation}

It is easy to see that $h\in A(B_X, Y)$ and $\|h(x)\|\le \|f\|+
\epsilon=\|h(\alpha x_0)\|$ for all $x\in B_X$. We claim that $h$ is
a strong peak function at $\alpha x_0$. Suppose that
$\lim_n\|h(x_n)\|=\|h\|$. For each $n$, we have
\[ \|h(x_n)\| \le \|f( x_n)\| + \epsilon |g( x^*(x_n))| \le \|f\| +
\epsilon=\|h\|.\] Hence $\lim_n\|f( x_n)\|=\|f\|$ and $\lim_n|g(
x^*(x_n))|=1$. Since $g$ is a peak function at $\alpha$,
$\{x^*(x_n)\}$ converges to $\alpha$. Now for any subsequence of
$\{x_n\}$, there is a further subsequence $\{y_k\}$ which converges
to $\eta x_0$ for some unit complex number $\eta$, because $\lim_n
\|f( x_n)\|= \|f\|$ and $f$ strongly attains its norm. Thus
$\alpha=\lim_k x^*(y_k) = \eta$. This means that every subsequence
of $\{x_n\}$ has a further subsequence converging to $\alpha x_0$,
so $\lim_n x_n = \alpha x_0$. Take $\tilde{g}(x) = \epsilon
g(x^*(x))\frac{f(\alpha x_0)}{\|f(\alpha x_0)\|}$. Then $\|\tilde
g\|\le \epsilon$ and $f+\tilde{g}$ is a strong peak function. Hence
we can conclude from Theorem~ \ref{thm:radonstrongnormattaining}
that the set of all strong peak functions in $A(B_X, Y)$ is dense in
$A(B_X, Y)$. The rest of proof follows from
Corollary ~\ref{cor:vectorGateaux}~(ii).\\

(ii) The proof follows from (i), Proposition~
\ref{prop:strongpeakfunctiondense} and the fact that every peak
point for $A(B_X, Y)$ is a complex extreme point of $B_X$ (see
\cite{Glo}).
\end{proof}

\begin{remark}\label{rem:peakfunction}
Notice that for any natural number $m$ the peak function $g$ at
$\alpha$ in \eqref{eq:peak1} can be chosen to be a polynomial
$g(\gamma) = (\bar\alpha\gamma + 1)^m/2^m $ of degree $m$. In
particular, the function $\tilde{g}(x) = \epsilon [\frac{\bar\alpha
(x^*(x)) + 1 }{2}]^m \frac{f(\alpha x_0)}{\|f(\alpha x_0)\|}$ is a
polynomial of degree $m$ and of rank $1$.
\end{remark}

Recall that a Banach space $X$ is said to be {\it locally uniformly
convex} if $x\in S_X$ and there is a sequence $\{x_n\}$ in $B_X$
satisfying $\lim_n \|x_n + x\|=2$, then $\lim_n \|x_n - x\|=0$.

Let $X$ be a complex Banach space. A point $x\in S_X$ is called a
{\it strong complex extreme point} of $B_X$ if for every
$\epsilon>0$, there exists $\delta>0$ such that
$$\sup_{0\le \theta \le 2\pi} \| x + e^{i\theta}y\| \ge 1 + \delta$$
for all $\|y\|\ge \epsilon$.  A complex Banach space $X$ is said to
be {\it locally uniformly $c$-convex} if every $x\in S_X$ is a
strong complex extreme point of $B_X$. Notice that if a complex
Banach space $X$ is locally uniformly convex, then $X$ is locally
uniformly $c$-convex. For more details on the local uniform
$c$-convexity, see \cite{DHM,Lee3}. A complex Banach space $X$ is
uniformly $c$-convex if for each $\epsilon>0$,
\[ H_X(\epsilon) = \inf\left\{ \sup_{0\le\theta\le 2\pi} \|x +
e^{i\theta}y \| -1\  :\  x\in S_X, y\in X, \|y\|\ge
\epsilon\right\}\] is strictly positive. It is easy to see that
every uniformly $c$-convex Banach space is locally uniformly
$c$-convex.

\begin{remark}
Let $X$ be a complex Banach space and let
\begin{align*}
A_{wu}(B_X)&=\{f\in A_u(B_X): f \text{ is weakly uniformly
continuous
on } B_X\}\\
A_{wb}(B_X)&=\{ f\in A_b(B_X): f\text{ is weakly continuous on }
B_X\}.
\end{align*} We
shall denote by $A_w(B_X)$ one of $A_{wb}(B_X)$ and $A_{wu}(B_X)$.
The proof of Theorem~\ref{thm:densitypeakRNP} and
Remark~\ref{rem:peakfunction} show that the set of all strong peak
functions for $A_{w}(B_X)$  is dense in $A_{w}(B_X)$ if $X$ has the
Radon-Nikod\'ym property.

It is a natural question that the set of all strong peak functions
in either $A(B_X)$ or $A_{w}(B_X)$ is dense, if $X$ has the {\it
analytic Radon-Nikod\'ym property}. The answer is negative in
$A_w(B_X)$ as observed in \cite{GLM}. Recall  that a complex Banach
space $X$ is said to have the {\it analytic Radon-Nikod\'ym
property} if for every bounded analytic function $f$ from the open
unit disc of $\mathbb{C}$ into $X$, it has the a.e. radial limits
\[ f(e^{i\theta}) = \lim_{r\uparrow 1} f(re^{i\theta}) \ \ \ \ a.e.\ \  \theta.\]
For more details on the analytic Radon-Nikod\'ym property, see
\cite{BD, GLM}.

Notice that $L_1[0,1]$ is uniformly $c$-convex and has the analytic
Radon-Nikod\'ym property (cf. \cite{Lee2}). Let $X=L_1[0,1]$. We
shall show that $A_w(B_X)$ does not contain any strong peak
function. Indeed, suppose that $f\in A_w(B_X)$ is a strong peak
function at $x$. For each $n\ge 1$, let
\[ U_n = \{ y\in B_X: |f(y)|> \|f\|-1/n\}.\] Then $U_n$ is a relative weak
neighborhood of $x$ for every $n$. Since $L_1[0,1]$ has the Daugavet
property, we can choose a sequence $\{x_n\}$ (see \cite{Sh}) such
that
\[ x_n \in U_n, \ \ \ \ \ \ \ \|x_n - x\|\ge 1, \ \ \ \ \ \ \ \forall n\ge 1.\]
This is a contradiction to that $f$ is a strong peak function at
$x$.
\end{remark}

\section{Density of norm-attaining elements in a subspace of $C_b(K,Y)$}

Let $X$ be a complex Banach space. An element $x\in B_X$ is said to
be a {\it strongly exposed point} for $B_X$ if there is a linear
functional $f\in B_{X^*}$ such that $f(x)=1$ and whenever there is a
sequence $\{x_n\}$ in $B_X$ satisfying $\lim_n \RE f(x_n)= 1$, we
get $\lim_n \|x_n -x\|= 0$. A set $\{x_{\alpha}\}$ of points on
$S_X$ is called {\it uniformly strongly exposed} (u.s.e.), if there
are a function $\delta(\epsilon)$ with $\delta(\epsilon)>0$ for
every $\epsilon
>0$, and a set $\{f_{\alpha}\}$ of elements of norm $1$ in $X^*$ such that
for every $\alpha$, $ f_{\alpha}(x_{\alpha}) = 1$, and for any $x$,
$$\|x\|\le 1 ~\text{and}~ \RE f_{\alpha}(x) \ge 1- \delta(\epsilon)
~\text{imply}~ \|x-x_{\alpha}\| \le \epsilon.$$ In this case we say
that $\{f_{\alpha}\}$ uniformly strongly exposes $\{x_{\alpha}\}$.
Lindenstrauss \cite[Proposition 1]{Li} showed that if $S_X$ is the
closed convex hull of a set of u.s.e. points, then $X$ has property
$A$, that is, for every Banach space $Y$, the set of norm-attaining
elements is dense in $L(X,Y)$, the Banach space of all bounded
operators of $X$ into $Y$. Modifying his argument and also applying
strong peak points instead of u.s.e. points, we study the density of
norm-attaining elements in a subspace of $C_b(K,Y)$. Notice that if
$S_X$ is the closed convex hull of a set $E$ of u.s.e. points, then
$E$ is a norming set for $L(X,Y)$.

\begin{thm}\label{thm:normattaining1}
Let $(K, d)$ be a complete metric space, $Y$ a Banach space and $A$
 a subspace of $C_b(K, Y)$. Assume that there exist a norming subset
$\{x_\alpha\}_\alpha\subset K$ for $A$ and a family
$\{\varphi_\alpha\}_\alpha$ of functions in $ C_b(K)$ such that each
$\varphi_\alpha$ is a  strong peak function at $x_\alpha$. Assume
also that $A$ contains $\varphi_\alpha^n \otimes y$ for each $y\in
Y$ and  $n\ge 1$. Then the set of norm-attaining elements of $A$ is
dense in $A$.
\end{thm}
\begin{proof}
We may assume that $\varphi_\alpha(x_\alpha)=1$ for each $\alpha$.
Let $f\in A$ with $\|f\|=1$ and $\epsilon$ with $0<\epsilon<1/3$ be
given. We choose a monotonically decreasing sequence
$\{\epsilon_k\}$ of positive numbers so that
\begin{equation}\label{eq1} 2\sum_{i=1}^\infty \epsilon_i <\epsilon,
\ \ \ 2\sum_{i=k+1}^\infty \epsilon_i < \epsilon_k^2,\ \ \
\epsilon_k < \frac 1{10k},\ \ \ k=1,2,\ldots\end{equation} We next
choose inductively sequences $\{f_k\}_{k=1}^\infty$,
$\{x_{\alpha_k}\}_{k=1}^\infty$ satisfying
\begin{align}
\label{eq2}&f_1=f\\
\label{eq3}&\|f_k(x_{\alpha_k})\|\ge \|f_k\|-\epsilon_k^2\\
\label{eq4}&f_{k+1}(x) = f_k(x) + \epsilon_k
\tilde{\varphi}_{\alpha_k}(x)\cdot f_k(x_{\alpha_k})\\
\label{eq9}|\tilde{\varphi}_{\alpha_k}(x)|&>1-1/k \ \ \text{ implies
} \ \ d(x,x_{\alpha_k})<1/k,
\end{align} where  $\tilde{\varphi}_{\alpha_j}$ is $\varphi^{n_j}_{\alpha_j}$ for some positive integer $n_j$.
Having chosen these sequences, we verify the following hold:
\begin{align}
\label{eq5}\|f_j - f_k\|&\le 2\sum_{i=j}^{k-1} \epsilon_i,\ \ \  \|f_k\|\le 4/3, &j<k,\ \ \  &k= 2, 3, \ldots\\
\label{eq6}\|f_{k+1}\|&\ge \|f_k\|+\epsilon_k\|f_k\| - 2\epsilon_k^2,&&k=1, 2,\ldots\\
\label{eq7}\|f_k\|&\ge \|f_j\|\ge 1, &j<k,\ \ \   &k=2, 3,\ldots\\
\label{eq8}|\tilde{\varphi}_{\alpha_j}(x_{\alpha_k})|&>1-1/j, &
j<k,\ \ \ &k=2, 3,\ldots.
\end{align}

Assertion (\ref{eq5}) is easy by using induction on $k$. By
(\ref{eq3}) and (\ref{eq4}), \begin{align*} \|f_{k+1}\|&\ge
\|f_{k+1} (x_{\alpha_k})\| = \|f_k(x_{\alpha_k})(1+\epsilon_k
\tilde{\varphi}_{\alpha_k}(x_{\alpha_k}))\| \\ &=
\|f_k(x_{\alpha_k})\|(1+\epsilon_k)\ge
(\|f_k\|-\epsilon_k^2)(1+\epsilon_k)\\&\ge
\|f_k\|+\epsilon_k\|f_k\|-2\epsilon_k^2,\end{align*} so the relation
(\ref{eq6}) is proved. Therefore (\ref{eq7}) is an immediate
consequence of (\ref{eq2}) and (\ref{eq6}). For $j<k$, by the
triangle inequality, (\ref{eq3}) and (\ref{eq5}),  we have
\begin{align*} \|f_{j+1}(x_{\alpha_k})\|&\ge
\|f_k(x_{\alpha_k})\|-\|f_k -f_{j+1}\|\\&\ge
\|f_k\|-\epsilon^2_k-2\sum_{i=j+1}^{k-1}\epsilon_j \ge
\|f_{j+1}\|-2\epsilon_j^2.\end{align*} Hence by (\ref{eq4}) and
(\ref{eq6}), \begin{align*}
\epsilon_j|\tilde{\varphi}_{\alpha_j}(x_{\alpha_k})| \cdot \|f_j\| +
\|f_j\|&\ge \|f_{j+1}(x_{\alpha_k})\|\ge \|f_{j+1}\|-2\epsilon^2_j\\
&\ge \|f_j\| + \epsilon_j \|f_j\|-4\epsilon_j^2,\end{align*} so that
\[ |\tilde{\varphi}_{\alpha_j}(x_{\alpha_k})| \ge 1-4\epsilon_j>1-1/j\] and
this proves (\ref{eq8}). Let $\hat{f}\in A$ be the limit of
$\{f_k\}$ in the norm topology. By (\ref{eq1}) and (\ref{eq5}),
$\|\hat{f}-f\|=\lim_n \|f_n-f_1\|\le2\sum_{i=1}^\infty \epsilon_i\le
\epsilon$ holds. The relations (\ref{eq9})and (\ref{eq8})  mean that
the sequence $\{x_{\alpha_k}\}$ converges to a point $\tilde{x}$,
say and  by (\ref{eq3}), we have $\|\hat{f}\|=\lim_n \|f_n\| =
\lim_n \|f_n(x_{\alpha_n})\| = \|\hat{f}(\tilde{x})\|$. Hence
$\hat{f}$ attains its norm. This concludes the proof.
\end{proof}

Let $A$ be the closed linear span of the constant 1 and $X^*$ as a
subspace of $C_b(B_X)$. Notice that if $X$ is locally uniformly
convex, then every element of $S_X$ is a strong peak point for $A$.
Therefore, every element of $S_X$ is a strong peak point for $A(B_X,
Y)$ for every complex Banach space $Y$, and $\rho A(B_X)$ is a
norming subset for $A(B_X,Y)$. Indeed, if $x\in S_X$, choose $x^*\in
S_{X^*}$ so that $x^*(x)=1$. Set $f(y) = \frac{x^*(y)+1}2$ for $y\in
B_X$. Then $f\in A$ and $f(x)=1$. If $\lim_n|f(x_n)|= 1$ for some
sequence $\{x_n\}$ in $B_X$, then $\lim_n x^*(x_n)=1$. Since
$|x^*(x_n) + x^*(x)|\le \|x_n + x\|\le 2$ for every $n$, $\|x_n +
x\|\to 2$ and $\|x_n -x\|\to 0$ as $n\to \infty$. Similarly it is
easy to see that every strongly exposed point for $B_X$ is a strong
peak point for $A$.

It was shown in \cite{CHL} that if a Banach sequence space $X$ is
locally uniformly $c$-convex and order continuous, then the set of
all strong peak points for $A(B_X)$ is dense in $S_X$. Therefore,
the set of all strong peak points for $A(B_X, Y)$ is dense in $S_X$
for every complex Banach space $Y$, and $\rho A(B_X)$ is a norming
subset for $A(B_X,Y)$. For the definition of a Banach sequence space
and order continuity, see \cite{CHL,FK,LT}. By the remarks above, we
get the following.

\begin{cor}\label{cor:analyticcase}
Suppose that $X$ and $Y$ are complex Banach spaces and $\rho A(B_X)$
is a norming subset for $A(B_X, Y)$. Then the set of norm-attaining
elements is dense in $A(B_X, Y)$. In particular, if $X$ is locally
uniformly convex, or if it is a locally uniformly $c$-convex, order
continuous Banach sequence space, then the set of norm-attaining
elements is dense in $A(B_X,Y)$.
\end{cor}

The complex Banach space $c_0$ renormed by Day's norm is locally
uniformly convex \cite{DGL,Day}, but it doesn't have the
Radon-Nikod\'{y}m property \cite{DU}. In addition, it is a locally
uniformly $c$-convex and order continuous Banach sequence space.

\begin{example}
A function $\varphi: \mathbb{R}\rightarrow [0,\infty]$ is said to be
an {\it Orlicz function} if $\varphi$ is even, convex continuous and
vanishing only at zero. Let $w=\{w(n)\}$ be a {\it weight sequence}
,that is, a non-increasing sequence of positive real numbers
satisfying $\sum_{n=1}^\infty w(n) =\infty$. Given a sequence $x$,
$x^*$ is the decreasing rearrangement of $|x|$.

An Orlicz-Lorentz sequence space $\lambda_{\varphi,w}$ consists of
all sequences $x=\{x(n)\}$ such that for some $\lambda>0$,
\[
\varrho_\varphi(\lambda x)= \sum_{n=1}^\infty \varphi(\lambda
x^*(n))w(n) < \infty,
\]
and equipped with the norm $\|x\| = \inf\{\lambda>0:
\varrho_\varphi(x/\lambda)\le 1\}$, which is a Banach sequence
space. We say an Orlicz function $\varphi$ satisfies {\it condition
$\delta_2$} $(\varphi\in \delta_2)$ if there exist $K>0$, $u_0>0$
such that $\varphi(u_0)>0$ and the inequality
\[ \varphi(2u) \le K \varphi(u)\] holds for $u\in [0, u_0]$.

If $\varphi\in \delta_2$, then $\lambda_{\varphi, w}$ is locally
uniformly $c$-convex \cite{CHL} and order continuous \cite{FK}.
Notice that if $\varphi(t) = |t|^p$ for $p\ge 1$ and $w=1$, then
$\lambda_{\varphi, w} = \ell_p$. The characterization of the local
uniform convexity of an Orlicz-Lorentz function space is given in
\cite{FK,HKM} and the characterization of the local uniform
$c$-convexity of a complex function space is given in \cite{Lee3}.
\end{example}

Extending the result of Lindenstrauss mentioned in the beginning of
this section, Pay\'{a} and Saleh \cite{PS} showed that if $B_X$ is
the closed absolutely convex hull of a set of u.s.e. points, then
the set of norm-attaining elements is dense in $L(^nX)$, the Banach
space of all bounded $n$-linear forms on $X$. We study a similar
question for the space of polynomials from $X$ into $Y$. In
particular, if a set of u.s.e. points on $S_X$ is a norming set for
the Banach space $P(^nX,Y)$ of all bounded $n$-homogeneous
polynomials from $X$ into $Y$, then the set of norm-attaining
elements is dense in $P(^nX,Y)$.

\begin{thm}Let $X$ and $Y$ be Banach spaces and $n\in \mathbb{N}$.
Suppose that a set $E$ of u.s.e. points on $S_X$ is a norming subset
of $P(^nX,Y)$. Then the set of all norm-attaining elements is dense
in $P(^nX,Y)$. Especially, if $E$ is dense in $S_X$, then the set of
norm-attaining elements is dense in $P(^nX,Y)$.

Moreover, if the set of strongly exposed points of $B_X$ is dense in
$S_X$, then the set of norm-attaining elements is dense in
$A(B_X,Y)$ for complex Banach spaces $X$ and $Y$.
\end{thm}
\begin{proof}Suppose that a set $E$ of u.s.e. points on $S_X$ is a norming subset
of $P(^nX,Y)$. Let $P\in P(^nX,Y)$, $\|P\|=1$,
 and $0 < \epsilon < 1/3 $ be
given. We first choose a monotonically decreasing sequence
$\{\epsilon_k\}$ of positive numbers so that
\begin{align}
\label{eq541}&4\sum_{i=1}^\infty \epsilon_i < \epsilon <
\frac{1}{3},\ \ \ \ 4\sum_{i=k+1}^\infty \epsilon_i <
\epsilon_{k}^2,\ \ \ \ \epsilon_k < \frac{1}{10k},\ \ \ \
k=1,2,\ldots.
\end{align}

Using induction, we next choose sequences $\{P_k \}_{k=1}^\infty$ in
$\mathcal{P}(^2X,Y)$, $\{x_k\}_{k=1}^\infty$ in $E$ and
$\{x^*_k\}_{k=1}^\infty$ in $S_{X^*}$ so that
\begin{align}
\label{eq542}&P_1=P\\
\label{eq543}&\|P_k(x_k)\|\ge \|P_k\|-\epsilon_k^2~  \text{  and }~
\|x_k\|=1,~~  x^*_k(x_k)=1,
\end{align}
where $\{x^*_k\}$  uniformly strongly
 exposes $\{x_k\}$,
\begin{align}
\label{eq544}&P_{k+1}(x)=P_k(x)+\epsilon_k\left(x^*_k(x)\right)^nP_k(x_k).
\end{align}

Having chosen these sequences, we see that the followings hold:
\begin{align}
\label{eq545}&\|P_j-P_k\|\le \frac{4}{3}\sum_{i=j}^{k-1}\epsilon_i,\
\ \ \ \ \ \
\|P_k\|\le \frac{4}{3},&j<k\\
\label{eq546}&\|P_{k+1}\| \ge \|P_k\|+\epsilon_k\|P_k\|
-\epsilon_k^2-\epsilon_k^3&\\
\label{eq547}&\|P_{k+1}\|
\le\|P_{k}\|+\epsilon_k|x_k^*(x_l)|^n\|P_k\|+\epsilon_k^2+2\cdot
\frac{4}{3}\sum_{i=k+1}^{l-1}\epsilon_i,& k+1 < l.
\end{align}

The assertion (\ref{eq545}) can easily be proved by induction and
(\ref{eq546}) follows directly from (\ref{eq544}). To see
(\ref{eq547}), for $k+1 <l$ we have
\begin{eqnarray*}
\|P_{k+1}\|&\le&  \|P_{l}\|+\|P_{k+1}-P_l\| \\
&\le&
\|P_l(x_l)\|+\epsilon_l^2+\frac{4}{3}\sum_{i=k+1}^{l-1}\epsilon_i\\
&\le&
\|P_{k}(x_l)\|+\epsilon_k|x_k^*(x_l)|^n\|P_k\|+\epsilon_k^2+2\cdot
\frac{4}{3}\sum_{i=k+1}^{l-1}\epsilon_i
\\
&\le& \|P_{k}\|+\epsilon_k|x_k^*(x_l)|^n\|P_k\|+\epsilon_k^2+2\cdot
\frac{4}{3}\sum_{i=k+1}^{l-1}\epsilon_i.
\end{eqnarray*}

By (\ref{eq545}), the sequence $\{P_k\}$ converges in the norm
topology to $Q\in P(^nX,Y)$ satisfying $\|P-Q\| < \epsilon$.

By (\ref{eq546}) and (\ref{eq547}) we have, for every $l> k+1$,
$$\epsilon_k\|P_k\|
-\epsilon_k^2-\epsilon_k^3\le
\epsilon_k|x_k^*(x_l)|^n\|P_k\|+2\epsilon_k^2,$$ and hence
$1-4\epsilon_k < |x_k^*(x_l)|^n.$

Since $A$ is uniformly strongly exposed,  $\{x_n\}$ has norm
convergent subsequence by Lemma~6 in \cite{A1}. Let $x_0$ be a limit
of that subsequence. Then we have $\|Q(x_0)\|=\|Q\|$. The rest of
the proof follows from Corollary~ \ref{cor:analyticcase} because
every strongly exposed point for $B_X$ is a strong peak point for
$A(B_X)$.
\end{proof}

Lindenstrauss \cite[Theorem 1]{Li} proved that the set of all
bounded linear operators of $X$ into $Y$ with norm-attaining second
adjoint is dense in $L(X,Y)$. In 1996 Acosta \cite{A} extended this
result to bilinear forms, and in 2002 Aron, Garcia and Maestre
\cite{ArGM} showed that this is also true for scalar-valued
$2$-homogeneous bounded polynomials. Recently, Acosta, Garcia and
Maestre \cite{AcGM} extended it to $n$-linear mappings.

We extend the result of \cite{ArGM} to the vector valued case by
modifying their proof, which is originally based on that of
Lindenstrauss. A bounded $n$-homogeneous polynomial $P\in P(^nX,Y)$
has an extension $\overline{P}\in P(^nX^{\ast\ast},Y^{\ast\ast})$ to
the bidual $X^{\ast\ast}$ of $X$, which is called the {\em
Aron-Berner extension} of $P$. In fact, $\overline{P}$ is defined in
the following way.

 Let
 $X_1, \cdots, X_n$ be  an arbitrary collection of  Banach spaces and
 let $\mathcal{L}(^n(X_1\times \cdots \times X_n))$ denote the space of bounded $n$-linear forms.
Given $z_i \in X^{**}_i$, $1\le i \le n$, define
 $\overline{z}_i$ from
 $\mathcal{L}(^n(X_1\times \cdots \times X_i\times X_{i+1}^{**}\times \cdots \times X_n^{**}))$ to
 $\mathcal{L}(^{n-1}(X_1\times \cdots \times X_{i-1}\times X_{i+1}^{**}\times \cdots \times X_n^{**}))$
by
\[\overline{z}_i(T)(x_1,\cdots,x_{i-1},x_{i+1}^{**},\cdots,x_n^{**})=\langle
z_i,T(x_1,\cdots,x_{i-1},\bullet,x_{i+1}^{**},\cdots,x_n^{**})\rangle,\]
where $T(x_1,\cdots,x_{i-1},\bullet,x_{i+1}^{**},\cdots,x_n^{**})$
is a linear functional on $X_i$ defined by $\bullet \mapsto
T(x_1,\cdots,x_{i-1},\bullet,x_{i+1}^{**},\cdots,x_n^{**})$  and
$\langle z, x^*\rangle$ is the duality between $X_{i}^{**}$ and
$X_{i}^*$. The map $\overline{z}_i$ is a bounded operator with norm
$\|z_i\|$. Now, given $T \in \mathcal{L}(^n(X_1\times \cdots \times
X_n ))$, define the extended $n$-linear form $\overline{T} \in
\mathcal{L}(^n(X_1^{**}\times \cdots \times X_n^{**} ))$ by
$$\overline{T}(z_1,\cdots z_n):= \overline{z}_1\circ \cdots \circ \overline{z}_n(T).$$

For a vector-valued $n$-linear mapping $L\in\mathcal{L}(^n(X_1\times
\cdots \times X_n),Y)$, define
$$\overline{L}(x_1^{**},\cdots,x_n^{**})(y^*)=\overline{y^*\circ L}(x_1^{**},\cdots,x_n^{**}),$$
where $x_i^{**}\in X_i^{**}$, $1\le i \le n$ and $y^* \in Y^*$. Then
$\overline{L}\in \mathcal{L}(^n(X_1^{**}\times \cdots \times
X_n^{**}),Y^{**})$ has the same norm as $L$. Let $S\in
\mathcal{L}_s(^nX,Y)$ be the symmetric $n$-linear mapping
corresponding to $P$, then $S$ can be extended to an $n$-linear
mapping $\overline{S} \in \mathcal{L}(^nX^{\ast\ast},Y^{\ast\ast})$
as described above. Then the restriction
$$\overline{P}(z) = \overline{S}(z,\ldots, z)$$
is called the Aron-Berner extension of $P$.  Given $z\in
X^{\ast\ast}$ and $w\in Y^\ast$, we have
$$\overline{P}(z)(w) = \overline{w\circ P}(z).$$
Actually this equality is often used as definition of the
vector-valued Aron-Berner extension based upon the scalar-valued
Aron-Berner extension.  Davie and Gamelin \cite[Theorem 8]{[DG]}
proved that $\|P\|=\|\overline{P}\|$. It is also worth to remark
that $\overline{S}$ is not symmetric in general.

\begin{thm}
Let $X$ and $Y$ be Banach spaces. The subset of $P(^2X,Y)$ each of
whose elements has the norm-attaining Aron-Berner extension is dense
in $P(^2X,Y)$.
\end{thm}
\begin{proof}Let $P\in \mathcal{P}(^2X,Y),~ \|P\|=1$, and
let $S$ be the symmetric bilinear mapping corresponding to $P$. Let
$\epsilon$ with $0 < \epsilon < 1/4 $ be given. We first choose a
monotonically decreasing sequence $\{\epsilon_k\}$ of positive
numbers which satisfies the following conditions:
\begin{align}
\label{551}&8\sum_{i=1}^\infty \epsilon_i < \epsilon < \frac{1}{4},
\ \ \ \  8\sum_{i=k+1}^\infty \epsilon_i < \epsilon_{k}^2  \ \ \
\text{ and }\ \ \  \epsilon_k < \frac{1}{10k},  \ \ \ \
k=1,2,\ldots.
\end{align}
Using induction, we next choose sequences $\{P_k \}_{k=1}^\infty$ in
$P(^2X,Y)$, $\{x_k\}_{k=1}^\infty$ in $S_X$ and
$\{f_k\}_{k=1}^\infty$ in $S_{Y^*}$ so that

\begin{align}
\label{552}&P_1=P, \ \ \ \  \|P\|=1\\
\label{553}&f_k(P_k(x_k))=\|P_k(x_k)\|\ge \|P_k\|-\epsilon_k^2\\
\label{554}&P_{k+1}(x)=P_k(x)+\epsilon_k\left(f_k(S_k(x_k,x))\right)^2P_k(x_k),
\end{align}
where each $S_k$ is the symmetric bilinear mapping corresponding to
$P_k$. Having chosen these sequences, we see that the following
hold:
\begin{align}
\label{555}&\|P_j-P_k\|\le
4\left(\frac{5}{4}\right)^3\sum_{i=j}^{k-1}\epsilon_i,~~~\ \ \ \ \
\|P_k\|\le
\frac{5}{4},~~~ & j<k\\
\label{556}&\|P_{k+1}\| \ge \|P_k\|+\epsilon_k\|P_k\|^3
-4\epsilon_k^2&\\
\label{557}&\|P_{j+1}(x_k)\|> \|P_{j+1}\|-2\epsilon_j^2,~~~  &j<k\\
\label{558}&|f_j(S_j(x_j,x_k))|^2\ge \|P_j\|^2-6\epsilon_j,~~~ &j<k
\end{align}

By (\ref{555}) and the polarization formula \cite{D}, the sequences
$\{P_k\}$ and $\{S_k\}$ converge in the norm topology to $Q$ and
$T$, say, respectively. Clearly $T$ is the symmetric bilinear
mapping corresponding to $Q$, and $\|P-Q\| < \epsilon.$

Let $\eta > 0$ be given. Then there exists $j_0 \in \mathbb{N}$ such
that
$$\|Q-P_j\| \le \|T-S_j\| < \eta~ \mbox{ for all } j \ge j_0,$$
hence $\|P_j\| \ge \|Q\| - \eta~ \mbox{ for all } j \ge j_0.$

By $$\|T-S_j\| \ge |f_j(T(x_j,x_k))-f_j(S_j(x_j,x_k))|$$ and
(\ref{558}), we have
\begin{eqnarray*}
|f_j(T(x_j,x_k))|& \ge &|f_j(S_j(x_j,x_k))| - \|T-S_j\|\\
   & \ge &\sqrt{\|P_j\|^2 - 6\epsilon_j} - \eta\\
   & \ge &\sqrt{(\|Q\|-\eta)^2 - 6\epsilon_j} - \eta
\end{eqnarray*}
for all $k > j \ge j_0$. Let $z \in X^{**}$ is a weak-$*$ limit
point of the sequence $\{x_k\}$. Then for all $j
 \ge j_0 $
$$\|\overline{T}(x_j,z)\| \ge \sqrt{(\|Q\|-\eta)^2-6\epsilon_j}-\eta.$$
Hence $\|\overline{T}(z,z)\| \ge \|Q\|-2\eta$. Since $\eta > 0$ is
arbitrary, we have
$$\|\overline{Q}(z)\|=\|\overline{T}(z,z)\|\ge \|Q\|=\|\overline{Q}\|.$$
\end{proof}

We finally investigate a version of Theorem~2 in \cite{Li} relating
with the complex convexity. Recall that a complex Banach space $X$
is said to be {\it strictly $c$-convex} if $S_X =
ext_\mathbb{C}(B_X)$.
\begin{thm}
Let $X$ be a Banach space with property $A$. Then
\begin{enumerate}
\item If $X$ is isomorphic to a strictly $c$-convex space, then
$B_X$ is the closed convex hull of its complex extreme points.

\item If $X$ is isomorphic to a locally uniformly $c$-convex space, then
$B_X$ is the closed convex hull of its strong complex extreme
points.
\end{enumerate}
\end{thm}
\begin{proof}
We prove only (2). We shall use the fact (\cite{Dil, DHM})  that
$x\in S_X$ is a strong complex extreme point of $B_X$ if and only if
for each $\epsilon>0$, there is $\delta>0$ such that
\[ \inf\left \{ \int_0^{2\pi} \|x + e^{i\theta}y\|^2 \ \frac{d\theta}{2\pi} \ :
\ y\in X, \|y\|\ge \epsilon \right\} \ge 1+\delta.\] For the proof
of (1), use the fact (\cite{Dil, DHM}) that $x\in S_X$ is a complex
extreme point of $B_X$ if and only if for any nonzero $y\in X$,
$\int_0^{2\pi} \|x + e^{i\theta}y\|^2 \ \frac{d\theta}{2\pi} > 1.$

 Let $C$ be the closed convex hull of the strong complex
extreme points of $B_X$. Suppose that $C\neq B_X$. Then there are
$f\in X^*$ with $\|f\|=1$ and $\delta,~0<\delta < 1$ such that
$|f(x)|<1-\delta$ for $x\in C$. Let $\||\cdot \||$ be a locally
uniformly $c$-convex norm on $X$, which is equivalent to the given
norm $\|\cdot \|$, such that $\||x\||\le \|x\|$ for $x\in X$. Let
$Y$ be the space $X\oplus_2 \mathbb{C}$ with the norm $\|(x,c)\|=
(\||x\||^2 + |c|^2)^{1/2}$. Then $Y$ is locally uniformly
$c$-convex. Otherwise, there exist $(x, c)\in S_{X\oplus_2
\mathbb{C}}$, $\epsilon>0$ and a sequence $\{(x_n, c_n)\}$ such that
for every $n\ge 1$, $\|(x_n, c_n)\|\ge \epsilon$ and
\[\lim_n \int_0^{2\pi} \|(x,c)  + e^{i\theta}(x_n, c_n)\|^2 \
\frac{d\theta}{2\pi}= 1.\] Since the norm is plurisubharmonic,
\begin{align*} 1 = \||x\||^2 + |c|^2 &\le \int_0^{2\pi} \|(x,c)  +
e^{i\theta}(x_n, c_n)\|^2\ \frac{d\theta}{2\pi}
\\&= \int_0^{2\pi} \||x+e^{i\theta} x_n\||^2\  \frac{d\theta}{2\pi}+
\int_0^{2\pi} |c+ e^{i\theta} c_n|^2 \ \frac{d\theta}{2\pi} \to
1.\end{align*} So
\[ \lim_{n\to\infty}\int_0^{2\pi} \||x+e^{i\theta} x_n\||^2 \ \frac{d\theta}{2\pi}
= \||x\||^2 \ \ \ \text{and} \ \ \ \lim_{n\to\infty}\int_0^{2\pi}
|c+ e^{i\theta} c_n|^2\  \frac{d\theta}{2\pi}=|c|^2.\] Since both
$(X, \||\cdot\||)$ and $\mathbb{C}$ are locally uniformly
$c$-convex, we get $\lim_n \||x_n\|| = \lim |c_n| =0$, which is a
contradiction to $\inf_n\|(x_n, c_n)\|\ge \epsilon$.

Let $V$ be the operator from $X$ into $Y$ defined by $Vx=(x,
Mf(x))$, where $M>2/\delta$. Then $V$ is an isomorphism (into) and
the same is true for every operator sufficiently close to $V$. We
have
\[ \|V\|\ge M, \ \ \ \ \|Vx\|\le (1+ (M-2)^2)^{1/2}\ \ \ \text{for } x\in C.\]
It follows that operators sufficiently close to $V$ cannot attain
their norm at a point belonging to $C$. To conclude its proof we
have only to show that if $T$ is an isomorphism (into) which attains
its norm at a point $x$ and if the range of $T$ is locally uniformly
$c$-convex, then $x$ is a strong complex extreme point of $B_X$.

We may assume that $\|Tx\|=\|T\|=1$. If $x$ is not a strong complex
extreme point, then there are $\epsilon>0$ and a sequence
$\{y_n\}\subset X$ such that $\|y_n\|\ge \epsilon$ for every $n$ and
\[\lim_n \int_0^{2\pi} \|x + e^{i\theta}y_n\|^2 \ \frac{d\theta}{2\pi}=1.\]
Then
\[ 1\le  \int_0^{2\pi} \|Tx + e^{i\theta}Ty_n\|^2 \ \frac{d\theta}{2\pi}\le
 \int_0^{2\pi} \|x + e^{i\theta}y_n\|^2 \ \frac{d\theta}{2\pi}\] shows
that $\{Ty_n\}$ converges to $0$, because the range of $T$ is
locally uniformly $c$-convex. Therefore, $\{y_n\}$ converges to $0$,
which is a contradiction.
\end{proof}

\section{Applications to a numerical boundary}

Let $\Pi_X = \{ (x, x^*)~ : ~ \|x\|=\|x^*\|=1=x^*(x)\}\subset S_X
\times S_{X^*}.$ We denote by $\tau$ the product topology of the
space $B_X\times B_{X^*}$, where the topologies on $B_X$ and
$B_{X^*}$ are the norm topology of $X$ and the weak-$*$ topology of
$X^*$, respectively. It is easy to see that $\Pi_X$ is a
$\tau$-closed subset of $B_X\times B_{X^*}$. Let $\pi_1$ be the
projection from $\Pi_X$ onto $S_X$ defined by $\pi_1(x,x^*) = x$ for
every $(x, x^*)\in \Pi_X$. It is not difficult to see that $\pi_1$
is a closed map.

The {\it spatial numerical range} of $f\in C_b(B_X, X)$ is defined
by $$ W(f) = \{x^*f(x) : (x^*, x)\in \Pi_X\}, $$ and the {\it
numerical radius} of $f\in C_b(B_X, X)$ is defined by $ v(f) =
\sup\{|\lambda| : \lambda\in W(f)\}.$ For a subspace $A\subset
C_b(B_X, X)$ we say that $B\subset \Pi_X$ is a numerical boundary
for $A$ if $$v(f) = \sup_{(x,x^*)\in B} |x^*(f(x))|,~~~\forall f\in
A,$$ and that $A$ has the numerical Shilov boundary if there is a
smallest closed numerical boundary for $A$. The numerical boundary
was introduced and studied in \cite{AK} for various Banach spaces,
and it was observed that the numerical Shilov boundary doesn't exist
for some Banach spaces. We first show that there exist the numerical
Shilov boundaries for most subspaces of $C_b(B_X, X)$ if $X$ is
finite dimensional. Notice that as a topological subspace of
$B_X\times B_{X^*}$, $\Pi_X$ is a compact metrizable space if $X$ is
finite dimensional.

\begin{thm} Let $X$ be a finite dimensional Banach space. Suppose
that a subspace $H$ of $C_b(B_X, X)$ contains the functions of the
form  \begin{equation}\label{eq:functionform} 1 \otimes x,\ \ \ \ \
y^*\otimes z, \ \ \ \ \forall x\in X,\ \ \forall z\in X,\ \ \
\forall y^*\in X^*.\end{equation} Then $H$ has the numerical Shilov
boundary.
\end{thm}
\begin{proof}
Consider the linear map $f\mapsto \tilde f$ from $H$ into $C(\Pi_X)$
defined by \[ \tilde f(x, x^*) = x^*f(x).\] Notice that $v(f) =
\|\tilde f\|$ for every $f\in H$. Let $\overline H$ be the closure
of the image $\tilde H$ in $C(\Pi_X)$. Then $\overline H$ is a
separable subspace of $C(\Pi_X)$.

We claim that $\overline H$ is separating. Let $(s, s^*) \neq (t,
t^*)\in \Pi_X$ and let $\alpha, \beta\in S_\mathbb{C}$. If $\alpha
t^* \neq \beta s^*$, then choose $x\in S_X$ such that $\alpha t^*(x)
\neq \beta s^*(x)$. Set $f= 1\otimes x \in H$. Then
\[\alpha \delta_{(t,t^*)}(\tilde f) =\alpha \tilde f (t,t^*) =
\alpha t^*(x) \neq \beta s^*(x) = \beta \tilde f(s,
s^*)=\beta\delta_{(s,s^*)}(\tilde f).\] If $\alpha t^* = \beta s^*$,
then $t\neq s$, and choose $z^*\in S_{X^*}$ such that $z^*(t) \neq
z^*(s)$. Set $f = z^* \otimes t \in H$. Then $\beta s^*(t) = \alpha
\neq 0$ and
\[ \alpha \tilde f (t, t^*) = \alpha z^*(t)t^*(t) =\beta z^*(t)s^*(t)\neq  \beta z^*(s) s^*(t)=
\beta \tilde f(s, s^*),\] hence $\alpha \delta_{(t,t^*)}(\tilde
f)\neq \beta\delta_{(s,s^*)}(\tilde f)$. Therefore $\overline H$ is
a separating separable subspace of $C(\Pi_X)$. By
Theorem~\ref{thm:Bishop}, there is the Shilov boundary $\partial
{\overline H}\subset \Pi_X$ for $\overline H$. It is clear that for
every $f\in H$,
\[ v(f) =\|\tilde f\|= \max_{(t, t^*)\in \partial {\overline H}} |t^*f(t)| .\]

We shall show that if $T\subset \Pi_X$ is a closed numerical
boundary for $H$, then $T$ is a closed boundary for $\overline H$.
Fix $g\in \overline{H}$ and choose a sequence $\{f_n\}_{n=1}^\infty$
in $H$ such that $\lim_n \|g-\tilde f_n\|= 0$. For each $n$, there
exists $(t_n, t_n^*)\in T$ such that $|t^*_nf_n (t_n)|= v(f_n)=
\|\tilde f_n\|$. So $\|g\|= \lim_n \|\tilde f_n\| = \lim_n |t^*_nf_n
(t_n)|$ and
\[|\tilde f_n(t_n, t_n^*) - g(t_n, t^*_n) |\le \|\tilde f_n - g\|\to
0.\] This shows that $\|g\|= \lim_n |g(t_n, t_n^*)|$ and
\[ \|g\| = \sup_{(t,t^*)\in T} |g(t, t^*)| =\max_{(t,t^*)\in T}
|g(t,t^*)|.\] Therefore, $T$ is a closed boundary for $\overline H$
and so $\partial {\overline H}$ is contained in $T$, which means
that $\partial {\overline H}$ is the smallest closed subset
satisfying
\[ v(f) = \max_{(t, t^*)\in \partial H} |t^*f(t)|, \ \ \ \ \forall
f\in H.\] The proof is done.
\end{proof}

\begin{example}
Let $X=\ell_\infty^2$ be the 2-dimensional space $\mathbb{C}^2$ with
the sup norm. Let $H$ be the subspace of $C_b(B_X, X)$ spanned by
all $f\otimes x$, $f\in X^*$ and $x\in X$. In fact, $H$ is
isometrically isomorphic to the Banach space $L(X)$ of bounded
linear operators from $X$ into $X$. It is easy to see that $v(T) =
\|T\|$ for $T\in H=L(X)$. Take
\begin{align*}S_1&= \{(x, x^*)~ :~  x=(x_1, 1)\in S_X,~|x_1| =1,~ x^*= (0, 1) \text{ or
} (\bar x_1, 0)] \},\\
S_2&=\{ (x, x^*)~ :~  x= (x_1, -1)\in S_X,~|x_1| =1,~ x^*=(0,-1)
\text{ or } (\bar x_1, 0)]\}.\end{align*} It is easy to see that for
each $T\in H$,
\[ v(T)= \|T\| = \sup_{(x, x^*)\in S_1} |x^*Tx| = \sup_{(x, x^*)\in S_2}
|x^*Tx|.\] However, $S_1$ and $S_2$ are disjoint closed subsets of
$\Pi_X$, so $H$ doesn't have the numerical Shilov boundary. In
particular, we cannot weaken the assumption of
Theorem~\ref{eq:functionform}.
\end{example}

Applying the Mazur theorem, we next prove the existence of the
numerical Shilov boundary for some subspaces of $C_b(B_X, X)$, when
$X$ is  separable.

\begin{thm}\label{thm:densitypeakpointnumericalbdry}
Let $X$ be a separable Banach space. Suppose that $A$ is a subspace
of $C_b(B_X)$ such that every element in $A$ is uniformly continuous
on $S_X$ and the set of all strong peak points for $A$ is dense in
$S_X$. If a subspace $H$ of $C_b(B_X, X)$ contains the functions of
the form: \begin{equation}\label{eq:span2} f\otimes y,\ \ \ \ \
\forall f\in A,\ \ \ \forall y\in X,\end{equation} then $H$ has the
numerical Shilov boundary. In particular, it is the set
\[ \overline{\{ (x, x^*): x \text{ is a smooth point of }
B_X\}}^{\ \tau} .\]
\end{thm}
\begin{proof}
Let $\Gamma = \{ (x, x^*): x \text{ is a smooth point of } B_X\}$.
We shall show that $\overline{\Gamma}^{\, \tau}$ is the  numerical
Shilov boundary for $H$. Notice that by Mazur's theorem, the set of
smooth points of $B_X$ is dense in $S_X$. Therefore, $\pi_1(\Gamma)$
is dense in $S_X$. By \cite[Theorem~2.5]{Pal},
$\overline{\Gamma}^{\, \tau}$ is a closed numerical boundary for
$H$, that is,
\[ v(f) = \max_{(t, t^*)\in \overline{\Gamma}^{\,\tau}} |t^*f(t)|, \ \ \ \ \forall
f\in H.\] Suppose that $C$ is a closed numerical boundary for $H$.
Then it is easy to see that $\pi_1(C)$ is a closed subset of $S_X$,
and $\pi_1(C)$ contains all strong peak points for $A$. Since the
set of all strong peak point for $A$ is dense in $S_X$,
$\pi_1(C)=S_X$. Therefore $\Gamma\subset C$, and hence
$\overline{\Gamma}^{\, \tau}\subset C$. This completes the proof.
\end{proof}

If $X$ is a smooth Banach space in Theorem
\ref{thm:densitypeakpointnumericalbdry}, then it is easily seen that
the numerical Shilov boundary for $H$ is $\Pi_X$, which is proved in
\cite{AK}.

\begin{cor}
Let a separable Banach space $X$ be locally uniformly convex and $A$
be a closed linear span of the constant 1 and $X^*$ as a subspace of
$C_b(B_X)$. Suppose that $H$ is a Banach space of uniformly
continuous functions from $B_X$ into $X$, which contains the
functions of the form \eqref{eq:span2}. Then the numerical Shilov
boundary for $H$ exists.
\end{cor}

\begin{cor}
Suppose that a Banach sequence space $X$ is locally uniformly
$c$-convex and also order continuous. If $H$ is a Banach space of
uniformly continuous functions from $B_X$ into $X$ which contains
the functions of the form:
\[ f\otimes y,\ \ \ \ \ \forall f\in A_u(B_X),\ \ \  \forall y\in X,\]
then $H$ has the numerical Shilov boundary.
\end{cor}

\bibliographystyle{amsplain}

\end{document}